\documentclass[a4paper,10pt,twoside]{article}
\usepackage{amssymb,amsmath,latexsym,geometry,fancyhdr,lineno,hyperref,titletoc,titlesec}
\contentsmargin{0pt}

\dottedcontents{section}[2em]{\vspace{-1mm}\small}{2em}{0pt}
\dottedcontents{subsection}[5em]{\vspace{-1mm}\small}{3em}{5pt}

\hypersetup{colorlinks,linkcolor= blue,citecolor=blue}

\newtheorem{theorem}{Theorem}[section]

\newtheorem{lemma}[theorem]{Lemma}

\newtheorem{definition}[theorem]{Definition}
\newtheorem{remark}{Remark}[section]

\def\proof{\mbox {\textbf{Proof.}~~}}
\numberwithin{equation}{section}
\allowdisplaybreaks[4]
\newcommand{\R}{\mathbb{R}}

\begin{document}
\title{{\bf\Large A mass supercritical problem revisited}}
\author{\\
{ \textbf{\normalsize Louis Jeanjean}}\\
{\it\small Laboratoire de Math\'{e}matiques (CNRS UMR 6623),}\\
{\it\small Universit\'{e} de Bourgogne Franche-Comt\'{e},}\\
{\it\small Besan\c{c}on 25030, France}\\
{\it\small e-mail: louis.jeanjean@univ-fcomte.fr}\\
\\
{ \textbf{\normalsize Sheng-Sen Lu}}\\
{\it\small Center for Applied Mathematics, Tianjin University,}\\
{\it\small Tianjin 300072, PR China}\\
{\it\small e-mail: sslu@tju.edu.cn}}
\date{}
\maketitle

{\bf\normalsize Abstract.} {\small
In any dimension $N\geq1$ and for given mass $m>0$, we revisit the nonlinear scalar field equation with an $L^2$ constraint:
    \begin{equation*}\tag{$P_m$}\label{eq:P_m}
      \left\{
             \begin{aligned}
               -\Delta u&=f(u)-\mu u\quad\text{in}~\mathbb{R}^N,\\
               \|u\|^2_{L^2(\mathbb{R}^N)}&=m,\\
               u&\in H^1(\mathbb{R}^N),
             \end{aligned}
      \right.
    \end{equation*}
where $\mu\in\mathbb{R}$ will arise as a Lagrange multiplier. Assuming only that the nonlinearity $f$ is continuous and satisfies weak mass supercritical conditions, we show the existence of ground states to \eqref{eq:P_m} and reveal the basic behavior of the ground state energy $E_m$ as $m>0$ varies. In particular, to overcome the compactness issue when looking for ground states, we develop robust arguments which we believe will allow treating other $L^2$ constrained problems in general mass supercritical settings. Under the same assumptions, we also obtain infinitely many radial solutions for any $N\geq2$ and establish the existence and multiplicity of nonradial sign-changing solutions when $N\geq4$. Finally we propose two open problems.
}




\pagestyle{fancy}
\fancyhead{} 
\fancyfoot{} 
\renewcommand{\headrulewidth}{0pt}
\renewcommand{\footrulewidth}{0pt}
\fancyhead[CE]{ \textsc{L. Jeanjean, S.-S. Lu}}
\fancyhead[CO]{ \textsc{a mass supercritical problem revisited}}

\fancyfoot[C]{\thepage}


\section{Introduction}\label{sect:introduction}
We are concerned with the nonlinear scalar field equation with an $L^2$ constraint:
  \begin{linenomath*}
    \begin{equation*}\tag{$P_m$}\label{P_m}
      \left\{
             \begin{aligned}
               -\Delta u&=f(u)-\mu u\quad\text{in}~\mathbb{R}^N,\\
               \|u\|^2_{L^2(\mathbb{R}^N)}&=m,\\
               u&\in H^1(\mathbb{R}^N).
             \end{aligned}
      \right.
    \end{equation*}
  \end{linenomath*}
Here $N\geq1$, $f\in C(\mathbb{R},\mathbb{R})$, $m>0$ is a given constant and $\mu\in\mathbb{R}$ will arise as a Lagrange multiplier. In particular $\mu\in\mathbb{R}$ does depend on the solution $u \in H^1(\mathbb{R}^N)$ and is not a priori given.

The main feature of \eqref{P_m} is that the desired solutions have an a priori prescribed $L^2$-norm. In the literature, solutions of this type are often referred to as normalized solutions. A strong motivation to study \eqref{P_m} is that it arises naturally in the search of stationary waves of nonlinear Schr\"{o}dinger equations of the following form
\begin{linenomath*}
    \begin{equation}\label{eq:equation-evolution}
     i \psi_t + \Delta \psi + g(|\psi|^2)\psi =0, \qquad \psi : \mathbb{R}_+ \times \mathbb{R}^N \to \mathbb{C}.
    \end{equation}
  \end{linenomath*}
Here by stationary waves we mean solutions of \eqref{eq:equation-evolution} of the special form $\psi(t,x) = e^{i\mu t} u(x)$ with a constant $\mu \in \mathbb{R}$ and a time-independent real valued function $u \in H^1(\mathbb{R}^N)$. The research of such type of equations started roughly forty years ago \cite{Be83-1,Be83-2,Lions84-1,Lions84-2,St77} and it now lies at the root of several models directly linked with current applications, such as nonlinear optics, the theory of water waves. For these equations, finding solutions with a prescribed $L^2$-norm is particularly relevant since this quantity is preserved along the time evolution.

Under mild conditions on $f$, one can introduce the $C^1$ functional
  \begin{linenomath*}
    \begin{equation*}\label{eq:functional}
      I(u):=\frac{1}{2}\int_{\mathbb{R}^N}|\nabla u|^2dx-\int_{\mathbb{R}^N}F(u)dx
    \end{equation*}
  \end{linenomath*}
on $H^1(\mathbb{R}^N)$, where $F(t):=\int^t_0f(\tau)d\tau$ for $t\in\mathbb{R}$. For any $m>0$, let
  \begin{linenomath*}
    \begin{equation*}
      S_m:=\Big\{u\in H^1(\mathbb{R}^N)~|~\|u\|^2_{L^2(\mathbb{R}^N)}=m\Big\}.
    \end{equation*}
  \end{linenomath*}
It is clear that solutions to \eqref{P_m} correspond to critical points of the functional $I$ constrained to the sphere $S_m$. Also, as may be well known, the study of \eqref{P_m} and the type of results one can expect depend on the behavior of the nonlinearity $f$ at infinity. In particular, this behavior determines whether $I$ is bounded from below on $S_m$ and so impacts on the choice of the approaches to search for constrained critical points.

In the present paper we shall focus on the mass supercritical case, that is, when $I$ is unbounded from below on $S_m$ for any $m>0$. Compared with the mass subcritical case, where the constrained functional $I_{|S_m}$ is bounded from below and coercive, more efforts are always needed in the study of the mass supercritical case. Indeed, even just aiming for an existence result, one has to identify a suspected critical level since it is no more possible to search for a global minimum of $I$ on $S_m$. Moreover, an arbitrary Palais-Smale sequence seems not necessarily bounded in $H^1(\mathbb{R}^N)$ let alone being strongly convergent up to a subsequence (and up to translations in $\mathbb{R}^N$ if necessary).

The first contribution to the mass supercritical case was made in \cite{Je97}. To make it more precise, we recall below the conditions introduced there.
\begin{itemize}
  \item[$(H0)$] $f : \mathbb{R} \to \mathbb{R}$ is continuous and odd.
  \item[$(H1)$] There exist $\alpha,\beta\in\mathbb{R}$ satisfying $2+4/N< \alpha \leq \beta < 2^*$ such that
                  \begin{linenomath*}
                    \begin{equation*}
                      0<\alpha F(t)\leq f(t)t\leq \beta F(t)\qquad\text{for any}~t\in\mathbb{R}\setminus\{0\},
                    \end{equation*}
                  \end{linenomath*}
                where $2^*:=\frac{2N}{N-2}$ for $N\geq3$ and $2^*:=+\infty$ when $N=1,2$.
  \item[$(H2)$] The function $\widetilde{F}(t):=f(t)t-2F(t)$ is of class $C^1$ and satisfies
                  \begin{linenomath*}
                    \begin{equation*}
                      \widetilde{F}'(t)t>\Big(2+\frac{4}{N}\Big)\widetilde{F}(t)\qquad\text{for any}~t\neq0.
                    \end{equation*}
                  \end{linenomath*}
\end{itemize}
In \cite{Je97}, under the conditions $(H0)$ and $(H1)$, the first author obtained a radial solution at a mountain pass value when $N\geq2$. Moreover, when $(H2)$ is also assumed, the existence of ground states was proved in any dimension $N\geq1$. Here by a ground state it is intended a solution $u$ to \eqref{P_m} that minimizes the functional $I$ among all the solutions to \eqref{P_m}:
  \begin{linenomath*}
    \begin{equation*}
      dI_{|S_m}(u)=0\qquad\text{and}\qquad I(u)=\inf \{ I(v)~|~dI_{|S_m}(v)=0 \}.
    \end{equation*}
  \end{linenomath*}
Afterwards, a multiplicity result was established by Bartsch and de Valeriola in \cite{BDV13}. When $N\geq2$ and $f$ satisfies $(H0)$ and $(H1)$, they derived infinitely many radial solutions from a fountain theorem type argument. In the very recent paper \cite{IT19}, Ikoma and Tanaka provided an alternative proof for this multiplicity result by exploiting an idea related to symmetric mountain pass theorems. One may also refer to \cite{BS17,BS18} for another proof which is based on a natural constraint approach but requires the additional assumption $(H2)$. More globally, the search of normalized solutions for problems that present a mass supercritical character is now a subject in full development.  We refer, for example, to \cite{AkWe19,BaJeSo,BS19,BaZhZo,BeBoJeVi17,BJL13,BCGL19,CiJe,NoTaVe19,So18,So19}.

Our aim in this work is to make a more in-depth study of \eqref{P_m} in the mass supercritical case. First, we relax some of the classical growth assumptions on $f$. For example, as one may observe, the condition $(H1)$ was required in all the previous papers \cite{BDV13,BS17,BS18,IT19,Je97}. In particular, the first part of this condition, i.e.,
  \begin{linenomath*}
    \begin{equation}\label{eq:key_1}
      \text{there exists}~\alpha >2+\frac{4}{N}~\text{such that}~0<\alpha F(t)\leq f(t)t~\text{for any}~t\neq0,
    \end{equation}
  \end{linenomath*}
was used in a technical but essential way not only in showing that the under study problem is mass supercritical but also in obtaining bounded constrained Palais-Smale sequences. We shall show that, under a weak version of the monotonicity condition $(H2)$, one can actually replace \eqref{eq:key_1} by a weaker and more natural mass supercritical condition. As a consequence, we manage to extend the previous results on the existence of ground states and the multiplicity of radial solutions. Moreover, we address new issues, such as the monotonicity of the ground state energy as a function of $m>0$ or the existence of infinitely many nonradial sign-changing solutions. Finally,  we stress that  all our results are obtained only assuming that the nonlinearity $f$, as any function  built on $f$,  is continuous.  In contrast to the cases where more regularity is assumed on the data, several classical tools are not available anymore and this forces us to develop, on several occasions, more robust proofs.  See in particular,  Remarks \ref{remark:robust}, \ref{remark:simplification}  and \ref{remark:f7}, in that direction.

Before stating the main results of this paper, let us present our conditions on $f$.
\begin{itemize}
    \item[$(f0)$] $f: \mathbb{R} \to \mathbb{R}$ is continuous.
    \item[$(f1)$] $\lim_{t\to0}f(t)/|t|^{1+4/N}=0$.
    \item[$(f2)$] When $N\geq3$, $\lim_{t\to\infty}f(t)/|t|^{\frac{N+2}{N-2}}=0$.

                  When $N=2$, $\lim_{t\to\infty}f(t)/e^{\gamma t^2}=0$ for any $\gamma>0$.
    \item[$(f3)$] $\lim_{t\to\infty}F(t)/|t|^{2+4/N}=+\infty$.
    \item[$(f4)$] $t\mapsto \widetilde{F}(t)/|t|^{2+4/N}$ is strictly decreasing on $(-\infty,0)$ and strictly increasing on $(0,\infty)$.
    \item[$(f5)$] When $N\geq3$, $f(t)t<\frac{2N}{N-2}F(t)$ for all $t\in\mathbb{R}\setminus\{0\}$.
  \end{itemize}

The conditions $(f0)-(f3)$ are somehow standard. They show that \eqref{P_m} is Sobolev subcritical but mass supercritical. $(f4)$ is a weaker version of $(H2)$ and plays a crucial role in this paper. In particular, it is under this condition that we can use $(f3)$ instead of \eqref{eq:key_1}. The condition $(f5)$ is weaker than the second part of $(H1)$ and only needed, in some dimension, to ensure that the Lagrange multipliers are positive. This proves crucial in our approaches to guarantee that certain bounded Palais-Smale sequences are strongly convergent up to a subsequence (and up to translations if necessary). At some points, we shall also make use of the following condition
  \begin{itemize}
    \item[$(f6)$] when $N\geq3$, $\lim_{t\to0}f(t)t/|t|^{\frac{2N}{N-2}}=+\infty$.
  \end{itemize}

As an example of the nonlinearity that fulfills $(f0)-(f6)$, setting
  \begin{linenomath*}
    \begin{equation*}
      \alpha_N:=
        \left\{
              \begin{aligned}
                1,\qquad&\qquad\text{for}~N=1,2,\\
                \frac{8}{N(N-2)},&\qquad\text{for}~N\geq3,
              \end{aligned}
        \right.
    \end{equation*}
  \end{linenomath*}
we have the odd continuous function
  \begin{linenomath*}
    \begin{equation*}
      f(t):=\Big[\big(2+\frac{4}{N}\big)\ln\big(1+|t|^{\alpha_N}\big)+\frac{\alpha_N|t|^{\alpha_N}}{1+|t|^{\alpha_N}}\Big]|t|^{4/N}t
    \end{equation*}
  \end{linenomath*}
with the primitive function $F(t):=|t|^{2+4/N}\ln \big(1+|t|^{\alpha_N}\big)$. In particular, this function does not satisfy \eqref{eq:key_1}. Since $(H1)$ implies that
  \begin{linenomath*}
    \begin{equation*}
      \min\{|t|^\alpha,|t|^\beta\}F(1)\leq F(t) \leq f(t)t \leq \beta \max\{|t|^\alpha,|t|^\beta\}F(1)\qquad \text{for any}~t\in \mathbb{R},
    \end{equation*}
  \end{linenomath*}
one can also see that the conditions $(f0)-(f6)$ are weaker than the previous ones $(H0)-(H2)$.

We are now in the position to present our main results. The first one concerns the existence of ground states and it reads as follows.

\begin{theorem}\label{theorem:groundstate}
  Let $N\geq1$ and $f$ satisfy $(f0)-(f4)$.
	\begin{itemize}
	\item[$(i)$] If the condition $(f5)$ is satisfied, then \eqref{P_m} admits a ground state for any $m>0$.
	\item[$(ii)$] Assume that $f$ is odd and that $(f5)$ holds for $N\geq 5$. Then \eqref{P_m} admits a positive ground state for any $m>0$.
	\end{itemize}
  In both cases, for any ground state the associated Lagrange multiplier $\mu$ is positive.
\end{theorem}

\begin{remark}\label{remark:(f4)}
  \begin{itemize}
    \item[$(i)$] As we shall see, the condition $(f4)$ permits to reduce the search of a ground state to a minimization problem set on a submanifold of $S_m$.  This is a specific feature of mass supercritical problems that conditions of the type of $(H2)$ or $(f4)$ appear necessary to consider the existence of a ground state. At least, there are no results so far without imposing such conditions or related ones as, for example, in \cite{CiJe,So18,So19}.
    \item[$(ii)$]  Under the condition $(f4)$, it is likely that considering our ground states as stationary solutions of the associated evolution problem, in the sense of \eqref{eq:equation-evolution}, one could prove that these ground states are unstable by blow-up in finite time. We refer to the classical paper \cite{BeCa} in that direction, see also \cite{Le} for further developments.
  \end{itemize}
\end{remark}

\begin{remark}\label{remark:(AR)}
In view of the role that condition \eqref{eq:key_1} plays in the constrained mass supercritical problems, it is reminiscent of the classical Ambrosetti-Rabinowitz condition introduced in \cite{AR73} for unconstrained superlinear problems. Indeed, our idea of weakening \eqref{eq:key_1} under the monotonicity condition $(f4)$ is somehow inspired by the papers \cite{Je99,LWZ06,LW04} where the authors demonstrated that under a Nehari type condition it is possible to find solutions without using the Ambrosetti-Rabinowitz condition.
\end{remark}

\begin{remark}\label{remark:(Radial)}
As one will see in Remark \ref{remark:addendum}, under the assumptions of Theorem \ref{theorem:groundstate} $(ii)$, we can actually obtain a positive radially symmetric ground state after having proved the existence of ground states. But this does not mean that, to prove Theorem \ref{theorem:groundstate} $(ii)$, one can work directly in the subspace of radially symmetric functions where additional compactness is available when $N\geq2$.  In fact,  in related problems, quite often it is shown at some point that it is not restrictive to work with sequences of functions which are Schwartz symmetric, see for example \cite[Lemma 4.2]{BeGeVi18}, for an illustration of this strategy. However, this possibility relies particularly on, for the present problem, a monotonicity property of the function $[f(t)t-(2+4/N)F(t)]/t^2$  that seems not to be guaranteed  under the setting of Theorem \ref{theorem:groundstate} $(ii)$. For more details in this direction, we refer to Remark \ref{remark:f7}.
\end{remark}

Let us now explain the strategy for the proof of Theorem \ref{theorem:groundstate} and highlight some of the difficulties encountered. First, for given $m>0$, we identify the suspected ground state energy
  \begin{linenomath*}
    \begin{equation}\label{eq:key_2}
      E_m :=\inf_{u \in \mathcal{P}_m} I(u),
    \end{equation}
  \end{linenomath*}
where $\mathcal{P}_m$ is the Pohozaev manifold defined by
  \begin{linenomath*}
    \begin{equation*}
      \mathcal{P}_m:=\Big\{u\in S_m~\Big|~P(u):=\int_{\mathbb{R}^N}|\nabla u|^2dx-\frac{N}{2}\int_{\mathbb{R}^N}\widetilde{F}(u)dx=0\Big\}.
    \end{equation*}
  \end{linenomath*}
As will be shown in Lemma \ref{lemma:I_4}, $\mathcal{P}_m$ is nonempty and $E_m>0$. Since $\mathcal{P}_m$ contains all the possible critical points of $I$ restricted to $S_m$, our task is to show that $E_m$ is a critical level of $I_{|S_m}$.

A difficulty appears when we try to construct a bounded Palais-Smale sequence of $I_{|S_m}$ at the level $E_m$. Indeed, under our assumptions on $f$, the information that a Palais-Smale sequence $\{u_n\}\subset S_m$ satisfies $P(u_n)=o_n(1)$, seems no more sufficient to prove its boundedness. To overcome this problem, in Lemma \ref{lemma:I_5}, we show that there exists one Palais-Smale sequence at the level $E_m$ which satisfies exactly
  \begin{linenomath*}
    \begin{equation*}
      P(u_n)=0 \qquad \text{for any}~n\geq1.
    \end{equation*}
  \end{linenomath*}
Since $I$ is coercive on $\mathcal{P}_m$ by Lemma \ref{lemma:I_4}, the boundedness follows. As one will see, our proof of Lemma \ref{lemma:I_5} borrows some arguments from Bartsch and Soave \cite{BS18,BS19}. However, since $\widetilde{F}$ is not required to be of class $C^1$, we need to adapt their argument by making use of techniques due to Szulkin and Weth \cite{SW09,SW10}. To be more precise, for any $u\neq0$ and $s\in\mathbb{R}$, let $(s\star u)(x):=e^{Ns/2}u(e^sx)$ for almost everywhere $x \in \mathbb{R}^N$. In Lemma \ref{lemma:I_3}, we show that a number $s(u)\in\mathbb{R}$ exists uniquely such that $P(s(u)\star u)=0$ and it is continuous as a mapping of $u\neq0$. Then, inspired by \cite[Proposition 2.9]{SW09}, we prove the $C^1$ regularity for the free functional
  \begin{linenomath*}
    \begin{equation*}
      \Psi(u):=I(s(u)\star u)=\frac{1}{2}e^{2s(u)}\int_{\mathbb{R}^N}|\nabla u|^2dx-e^{-Ns(u)}\int_{\mathbb{R}^N}F(e^{Ns(u)/2}u)dx
    \end{equation*}
  \end{linenomath*}
on $H^1(\mathbb{R}^N)\setminus\{0\}$, see Lemma \ref{lemma:Psi}. After that, we manage to produce the desired Palais-Smale sequence by adapting the arguments of \cite[Proposition 3.9]{BS19} to the $C^1$ constrained functional $J:=\Psi_{|S_m}$, see Lemma \ref{lemma:J_1} and the proof of Lemma \ref{lemma:I_5}.

Since we search for solutions having a given $L^2$-norm, we must deal with a possible lack of compactness for the above bounded Palais-Smale sequence $\{u_n\} \subset \mathcal{P}_m$. It is not difficult to see that, up to a subsequence and up to translations in $\mathbb{R}^N$, the sequence $\{u_n\}$ has a nontrivial weak limit $u \in H^1(\mathbb{R}^N)$. Moreover, $s:=\|u\|^2_{L^2(\mathbb{R}^N)}\in (0,m]$,
  \begin{linenomath*}
    \begin{equation}\label{eq:key_3}
      -\Delta u = f(u)-\mu u \qquad\text{for some}~\mu\in\mathbb{R},
    \end{equation}
  \end{linenomath*}
and $u\in\mathcal{P}_s$. By a technical argument, we also prove that $\lim_{n\to\infty}I(u_n-u)\geq 0$ and hence
  \begin{linenomath*}
    \begin{equation*}
      E_m=\lim_{n\to\infty}I(u_n)=I(u)+\lim_{n\to\infty}I(u_n-u)\geq I(u)\geq E_s.
    \end{equation*}
  \end{linenomath*}
Clearly, the compactness would be proved if one can show that $E_m<E_s$ for any $s\in (0,m)$ or more globally that $E_m$ is strictly decreasing as a function of $m>0$. Therefore, from this observation, the study of the monotonicity of the function $m\mapsto E_m$ arises naturally as a fundamental problem. In this direction, we have the following theorem which also reveals some other basic properties of $E_m$.
\begin{theorem}\label{theorem:Em}
  Let $N\geq1$ and $f$ satisfy $(f0)-(f4)$. Then the function $m \mapsto E_m$ is positive, continuous, nonincreasing and $\lim_{m \to 0^+} E_m=+\infty$. Moreover, when $N=1,2$, we have
    \begin{itemize}
      \item[$(i)$] $E_m$ is strictly decreasing in $m>0$,
      \item[$(ii)$] $\lim_{m\to\infty} E_m=0$.
    \end{itemize}
  When $N \geq 3$, Items $(i)$ and $(ii)$ hold if $f$ also satisfies $(f5)$ and $(f6)$ respectively. In particular, when $f$ is odd and $N=3,4$, we have Item $(i)$ without assuming $(f5)$.
\end{theorem}

One should however note that the strict decrease of $E_m$ can be established only after having proved Theorem \ref{theorem:groundstate} and so it actually does not play a role in recovering the compactness. As we shall see, what really works in practice are the basic property that $E_m$ is nonincreasing and a companion result Lemma \ref{lemma:Em_strictlymonotonic}. They permit to reduce the problem of strong convergence to the one of showing that the Lagrange multiplier $\mu\in\mathbb{R}$ in \eqref{eq:key_3} is positive. For more details, we refer to the last part of the proof of Lemma \ref{lemma:decomposition}.

\begin{remark}\label{remark:Em}
  \begin{itemize}
    \item[$(i)$] To prove Theorem \ref{theorem:Em} (and Lemma \ref{lemma:Em_strictlymonotonic}), we develop robust arguments which we believe will allow treating other $L^2$ constrained problems in general mass supercritical settings. In this direction, we refer to Remark \ref{remark:robust} for more details.
    \item[$(ii)$] When $(f6)$ does not hold, the limit $\lim_{m\to\infty} E_m$ can be positive, see Remark \ref{remark:E_infinity}.
    \item[$(iii)$] From a later result, Lemma \ref{lemma:J}, and the below characterization
          \begin{linenomath*}
            \begin{equation*}
              E_m :=\inf_{u \in \mathcal{P}_m} I(u)=\inf_{u \in S_m}J(u),
            \end{equation*}
          \end{linenomath*}
        one can see that any minimizer $u \in \mathcal{P}_m$ of \eqref{eq:key_2} is a solution and thus a ground state to \eqref{P_m}. However, despite this fact, it seems not a good choice to prove Theorem \ref{theorem:groundstate} by solving directly the minimization problem. Indeed, when our Palais-Smale sequence $\{u_n\} \subset \mathcal{P}_m$ is replaced by an arbitrary minimizing sequence of \eqref{eq:key_2}, up to a subsequence and up to translations in $\mathbb{R}^N$, it still has a nontrivial weak limit $u \in H^1(\mathbb{R}^N)$ with $s:=\|u\|^2_{L^2(\mathbb{R}^N)}\in (0,m]$ but the information $u \in \mathcal{P}_s$ and $\lim_{n\to\infty}I(u_n-u)\geq0$ seems now out of reach. This also explains why we introduce Palais-Smale sequences to solve the minimization problem \eqref{eq:key_2}.
  \end{itemize}
\end{remark}

Our next result concerns the existence of infinitely many radial solutions when $N\geq2$.
\begin{theorem}\label{theorem:radial}
  Assume that $N\geq2$ and $f$ is odd satisfying $(f0)-(f5)$. Then \eqref{P_m} has infinitely many radial solutions $\{u_k\}^\infty_{k=1}$ for any $m>0$. In particular,
    \begin{linenomath*}
      \begin{equation*}
        I(u_{k+1})\geq I(u_k)>0\qquad\text{for each}~k\in\mathbb{N}^+
      \end{equation*}
    \end{linenomath*}
  and $I(u_k)\to+\infty$ as $k\to\infty$.
\end{theorem}
\begin{remark}
  It is clear that Theorem \ref{theorem:radial} extends \cite[Theorem 1.4]{BS17} where the stronger conditions $(H0)-(H2)$ were assumed. Despite the fact that the multiplicity result in \cite{BDV13,IT19} did not require a monotonicity condition like $(f4)$, our Theorem \ref{theorem:radial} is not a special case of theirs. Indeed, the conditions and thus the results are mutually non-inclusive and the methods are also different.
\end{remark}

We now give some ideas of the proof of Theorem \ref{theorem:radial}. We work in $H^1_r(\mathbb{R}^N)$, the space of radially symmetric functions in $H^1(\mathbb{R}^N)$. Since the constrained functional $J:=\Psi_{|S_m}$ is even, using the genus theory, it is not difficult to define an infinite sequence of minimax values $E_{m,k}$. In particular, $E_{m,k}$ is positive and nondecreasing in $k\geq1$, see Lemma \ref{lemma:GkEmk}. By a similar argument as in the proof of Theorem \ref{theorem:groundstate}, we establish Lemma \ref{lemma:J_2} which will be used to ensure the existence of a Palais-Smale sequence $\{u^k_n\}^\infty_{n=1} \subset \mathcal{P}_m \cap H^1_r(\mathbb{R}^N)$ for the constrained functional $I_{|S_m\cap H^1_r(\mathbb{R}^N)}$ at each level $E_{m,k}$. Regarding the problem of strong convergence, it can be solved by the compactness result Lemma \ref{lemma:compact} whose proof uses essentially the fact that the inclusion $H^1_r(\mathbb{R}^N)\hookrightarrow L^p(\mathbb{R}^N)$ is compact for any $2<p<2^*$. To conclude the proof, we also need to show that $E_{m,k}$ is unbounded. Since the Pohozaev manifold $\mathcal{P}_m$ is only a topological manifold, it seems no more possible to prove this by a standard genus type argument for $I_{|\mathcal{P}_m}$. Fortunately, inspired by \cite[Theorem 9]{Be83-2}, we manage to justify this key point by developing a new argument, see Lemma \ref{lemma:unbounded} and the proof of Lemma \ref{lemma:Emk}.

In the last part of our study, we are interested in the construction of nonradial sign-changing solutions to \eqref{P_m}. To state our results in this direction, we introduce some notations at first. Assume that $N\geq4$ and $2\leq M\leq N/2$. Let us fix a transformation $\omega\in \mathcal{O}(N)$ such that $\omega(x_1,x_2,x_3)=(x_2,x_1,x_3)$ for any $x_1,x_2\in\mathbb{R}^M$ and $x_3\in\mathbb{R}^{N-2M}$, where $x=(x_1,x_2,x_3)\in\mathbb{R}^N=\mathbb{R}^M\times\mathbb{R}^M\times\mathbb{R}^{N-2M}$. We define the Sobolev space of odd functions
  \begin{linenomath*}
    \begin{equation*}
      X_\omega:=\big\{u\in H^1(\mathbb{R}^N)~|~u(\omega x)=-u(x)~\text{for a.e.}~x\in\mathbb{R}^N\big\},
    \end{equation*}
  \end{linenomath*}
which clearly does not contain nontrivial radial functions. Let $H^1_{\mathcal{O}_1}(\mathbb{R}^N)$ denote the subspace of invariant functions with respect to $\mathcal{O}_1$, where $\mathcal{O}_1:=\mathcal{O}(M)\times\mathcal{O}(M)\times \text{id}\subset \mathcal{O}(N)$ acts isometrically on $H^1(\mathbb{R}^N)$. We also consider $\mathcal{O}_2:=\mathcal{O}(M)\times\mathcal{O}(M)\times\mathcal{O}(N-2M)\subset \mathcal{O}(N)$ acting isometrically on $H^1(\mathbb{R}^N)$ with the subspace of invariant functions denoted by $H^1_{\mathcal{O}_2}(\mathbb{R}^N)$. Here we agree that the components corresponding to $N-2M$ do not exist when $N=2M$. It is clear that $H^1_{\mathcal{O}_2}(\mathbb{R}^N)$ is in general a subspace of $H^1_{\mathcal{O}_1}(\mathbb{R}^N)$ but coincides with the latter when $N=2M$.

Now, for notational convenience, we set
  \begin{linenomath*}
    \begin{equation*}
      X_1:=H^1_{\mathcal{O}_1}(\mathbb{R}^N)\cap X_\omega \qquad \text{and} \qquad X_2:=H^1_{\mathcal{O}_2}(\mathbb{R}^N)\cap X_\omega.
    \end{equation*}
  \end{linenomath*}
In any dimension $N\geq4$, we have the following existence result of nonradial solutions.
\begin{theorem}\label{theorem:nonradial1}
  Assume that $N\geq4$ and $f$ is odd satisfying $(f0)-(f5)$.  Then \eqref{P_m} has a nonradial solution $v\in X_1$ for any $m>0$. In particular, $v$ changes signs, minimizes $I$ among all the solutions of \eqref{P_m} belonging to $X_1$ and $I(v)>2E_m$.
\end{theorem}
\begin{remark}
  The nonradial solution obtained here can be regarded as a ground state within the subspace $X_1$.  When $N\geq4$ and $N-2M\neq0$, $X_1$ does not embed compactly into any $L^p(\mathbb{R}^N)$ and so, in this case, finding that nonradial solution is similar to the search of a ground state in $H^1(\mathbb{R}^N)$.
\end{remark}

When $N=4$ or $N\geq 6$, we can choose $M\geq2$ such that $N-2M\neq 1$. In this case, we can obtain infinitely many nonradial solutions in $X_2$.
\begin{theorem}\label{theorem:nonradial2}
  Assume that $N=4$ or $N\geq6$, $N-2M\neq1$, and $f$ is odd satisfying $(f0)-(f5)$.  Then \eqref{P_m} possesses infinitely many nonradial solutions $\{v_k\}^\infty_{k=1}\subset X_2$ for any $m>0$. In particular, all these nonradial solutions change signs,
    \begin{linenomath*}
      \begin{equation*}
        I(v_{k+1})\geq I(v_k)>0\qquad\text{for each}~k\in\mathbb{N}^+
      \end{equation*}
    \end{linenomath*}
  and $I(v_k)\to+\infty$ as $k\to\infty$.
\end{theorem}

For the free nonlinear scalar field equation
\begin{linenomath*}
    \begin{equation*}
      - \Delta u = f(u) - \mu u, \qquad u \in H^1(\mathbb{R}^N),
    \end{equation*}
  \end{linenomath*}
the question of the existence of nonradial solutions was raised by Berestycki and Lions \cite[Section 10.8]{Be83-2} and it has been much studied over the past few decades. A positive answer was first given by Bartsch and Willem \cite{Ba93} in dimension $N=4$ and $N\geq6$. The idea of working within the subspaces as $X_2:=H^1_{\mathcal{O}_2}(\mathbb{R}^N)\cap X_\omega$ originates from \cite{Ba93}. Later on, Lorca and Ubilla \cite{Lo04} coped with the case $N=5$ by introducing the $\mathcal{O}_1$ action on $H^1(\mathbb{R}^N)$. In a more recent work \cite{MPW12}, using Lyapunov-Schmidt reduction methods, Musso, Pacard and Wei constructed nonradial solutions for any dimension $N\geq2$. We also would like to mention the most recent advance made in \cite{Me17}. Under the general Berestycki-Lions conditions, Mederski proved the existence and multiplicity of nonradial solutions when $N\geq4$; see also \cite{JL20} for an alternative proof with more elementary arguments.

However, in sharp contrast to the above free case, the study of nonradial solutions is almost unexplored in the literature for the constrained problem \eqref{P_m}. The first and currently the only paper to deal with normalized nonradial solutions is \cite{JL19} where the authors considered the mass subcritical case in a very general setting. In the present paper, regarding the issue of normalized nonradial sign-changing solutions, we somehow extend the existence and multiplicity results of \cite{JL19} to the mass supercritical case. To prove Theorems \ref{theorem:nonradial1} and \ref{theorem:nonradial2}, we shall adapt the arguments of Theorems \ref{theorem:groundstate} and \ref{theorem:radial}.

The remaining part of this paper is organized as follows. We present in Section \ref{sect:preliminaries} some preliminary results and then study in Section \ref{sect:Em} some basic properties of the function $m \mapsto E_m$. In Section \ref{sect:groundstate} we complete the proofs of Theorems \ref{theorem:groundstate} and \ref{theorem:Em}. Sections \ref{sect:radial} and \ref{sect:nonradial} deal with the existence of infinitely many radial solutions and the existence and multiplicity of nonradial sign-changing solutions respectively. Finally, in Section \ref{sect:final}, we justify Remark \ref{remark:Em} $(ii)$ and propose two open problems.

\begin{remark}\label{remark:addendum}
After the completion of this work, we were informed by J. Mederski of the manuscript \cite{BM20} which has some overlap with ours, at the level of Theorems \ref{theorem:groundstate} and \ref{theorem:Em}. In \cite{BM20}, for $N \geq 3$, a new view to the problem of the existence of a ground state is introduced.  The basic idea is to  transform the problem of the existence of a ground state on $S_m$ to the search of a global minima for $I$ on
\begin{linenomath*}
    \begin{equation*}
      \mathcal{P}:=\Big\{u\in H^1(\mathbb{R}^N)\setminus\{0\}~\big|~\|u\|^2_{L^2(\mathbb{R}^N)}\leq m~~\text{and}~~P(u)=0\Big\}.
    \end{equation*}
  \end{linenomath*}
This interesting approach relies, at least so far, on stronger regularity assumptions and in particular  the function $\widetilde{F}$ needs  to be of class $C^1$, see \cite[Theorem 1.1]{BM20}. The approach of \cite{BM20}  does not thus permit to recover the results of Theorems \ref{theorem:groundstate} and \ref{theorem:Em} in full generality. Nevertheless the paper \cite{BM20} proved  useful to show that, for odd functions satisfying $(f0)-(f4)$, if $E_m$ is reached then \eqref{P_m} admits a Schwarz symmetric ground state.  Indeed, let $v\in \mathcal{P}_m$ be a minimizer of $E_m$ and define $\tilde{v}:= |v|^*$ as the Schwarz symmetrization of $|v|$. It follows that $\tilde{v} \in S_m$, $I(\tilde{v}) \leq I(v) =E_m$ and $P(\tilde{v}) \leq P(v) = 0$. Clearly, if $P(\tilde{v})=0$, then $I(\tilde{v})=E_m$ and we complete the proof in view of Remark \ref{remark:Em} $(iii)$. To this end, we assume by contradiction that $P(\tilde{v}) <0$. Inspired by \cite[Lemma 2.7]{BM20} we observe that there exists $t:= t(\tilde{v})>1$ such that $P(\tilde{v}(t\cdot))=0$ and $m' := \|\tilde{v}(t\cdot)\|^2_{L^2(\mathbb{R}^N)} <m$. Since Theorem \ref{theorem:Em} implies that $E_m >0$ and that the function $m \mapsto E_m$ is nonincreasing, we have
\begin{linenomath*}
    \begin{equation*}
      \begin{split}
        E_m \leq E_{m'} \leq I(\tilde{v}(t\cdot))
        & = t^{-N} \left[\frac{N}{4}\int_{\R^N}\Big(\widetilde{F}(\tilde{v})- \frac{4}{N}F(\tilde{v})\Big)dx \right]\\
        &<  \frac{N}{4}\int_{\R^N}\Big(\widetilde{F}(\tilde{v})- \frac{4}{N}F(\tilde{v})\Big)dx\\
        & = \frac{N}{4}\int_{\R^N}\Big(\widetilde{F}(v)- \frac{4}{N}F(v)\Big)dx= I(v) = E_m.
      \end{split}
    \end{equation*}
  \end{linenomath*}
This contradiction proves that $P(\tilde{v})=0$ and thus $\tilde{v} \in S_m$ is a Schwarz symmetric ground state.
\end{remark}

\section{Preliminary results}\label{sect:preliminaries}
In this section we prepare several technical results for the proofs of our main Theorems \ref{theorem:groundstate}--\ref{theorem:nonradial2}. For notational convenience, we set
  \begin{linenomath*}
    \begin{equation*}
      B_m:=\Big\{u\in H^1(\mathbb{R}^N)~|~\|u\|^2_{L^2(\mathbb{R}^N)}\leq m\Big\}
    \end{equation*}
  \end{linenomath*}
for any $m>0$. The first technical result reads as follows and will be often used in the sequel.
\begin{lemma}\label{lemma:I_1}
  Assume that $N\geq1$ and $f$ satisfies $(f0)-(f2)$. Then the following statements hold.
    \begin{itemize}
      \item[$(i)$] For any $m>0$, there exists $\delta=\delta(N,m)>0$ small enough such that
                     \begin{linenomath*}
                       \begin{equation*}
                         \frac{1}{4}\int_{\mathbb{R}^N}|\nabla u|^2dx\leq I(u)\leq \int_{\mathbb{R}^N}|\nabla u|^2dx
                       \end{equation*}
                     \end{linenomath*}
                   for all $u\in B_m$ satisfying $\|\nabla u\|_{L^2(\mathbb{R}^N)}\leq \delta$.
      \item[$(ii)$] Let $\{u_n\}$ be a bounded sequence in $H^1(\mathbb{R}^N)$. If $\lim_{n\to\infty}\|u_n\|_{L^{2+4/N}(\mathbb{R}^N)}=0$, then
                     \begin{linenomath*}
                       \begin{equation*}
                         \lim_{n\to\infty}\int_{\mathbb{R}^N}F(u_n)dx=0=\lim_{n\to\infty}\int_{\mathbb{R}^N}\widetilde{F}(u_n)dx.
                       \end{equation*}
                     \end{linenomath*}
      \item[$(iii)$] Let $\{u_n\},\{v_n\}$ be bounded sequences in $H^1(\mathbb{R}^N)$ and $\lim_{n\to\infty}\|v_n\|_{L^{2+4/N}(\mathbb{R}^N)}=0$. Then
                     \begin{linenomath*}
                       \begin{equation}\label{eq:fuv}
                         \lim_{n\to\infty}\int_{\mathbb{R}^N}f(u_n)v_ndx=0.
                       \end{equation}
                     \end{linenomath*}
    \end{itemize}
\end{lemma}
\proof  We provide a full proof for Item $(i)$ but, for saving space, we only consider the case $N=2$ for the remaining two items.

$(i)$ We only need to show that there exists $\delta=\delta(N,m)>0$ small enough such that
  \begin{linenomath*}
    \begin{equation}\label{eq:F1}
      \int_{\mathbb{R}^N}|F(u)|dx\leq \frac{1}{4}\int_{\mathbb{R}^N}|\nabla u|^2dx\qquad\text{for any}~u\in B_m~\text{with}~\|\nabla u\|_{L^2(\mathbb{R}^N)}\leq \delta.
    \end{equation}
  \end{linenomath*}

We first prove \eqref{eq:F1} when $N\geq3$. Let $\varepsilon>0$ be arbitrary. By $(f0)-(f2)$, there exists $C_\varepsilon>0$ such that $|F(t)|\leq \varepsilon |t|^{2+\frac{4}{N}}+C_\varepsilon |t|^{\frac{2N}{N-2}}$ for all $t\in\mathbb{R}$. For any $u\in B_m$, using also Gagliardo-Nirenberg inequality, one then has
  \begin{linenomath*}
    \begin{equation*}
      \begin{split}
        \int_{\mathbb{R}^N}|F(u)|dx
        &\leq \varepsilon \int_{\mathbb{R}^N}|u|^{2+\frac{4}{N}}dx+C_\varepsilon\int_{\mathbb{R}^N}|u|^{\frac{2N}{N-2}}dx\\
        &\leq \varepsilon C_N m^{\frac{2}{N}}\int_{\mathbb{R}^N}|\nabla u|^2dx+C_\varepsilon C'_N\Big(\int_{\mathbb{R}^N}|\nabla u|^2dx\Big)^{\frac{N}{N-2}}\\
        &=\Big[\varepsilon C_N m^{\frac{2}{N}}+C_\varepsilon C'_N\Big(\int_{\mathbb{R}^N}|\nabla u|^2dx\Big)^{\frac{2}{N-2}}\Big]\int_{\mathbb{R}^N}|\nabla u|^2dx,
      \end{split}
    \end{equation*}
  \end{linenomath*}
where $C_N, C'_N>0$ depend only on $N$. Clearly, we obtain \eqref{eq:F1} by setting
  \begin{linenomath*}
    \begin{equation*}
      \varepsilon:=\frac{1}{8C_N m^{2/N}}\qquad\text{and}\qquad \delta:=\Big(\frac{1}{8C_\varepsilon C'_N}\Big)^\frac{N-2}{4}.
    \end{equation*}
  \end{linenomath*}

When $N=2$, for $\gamma :=1/(m+1)$ and any $\varepsilon>0$, by $(f0)-(f2)$, one can find $C'_\varepsilon>0$ such that $|f(t)|\leq \varepsilon |t|^3+C'_\varepsilon |t|^5e^{\gamma t^2/2}$ for all $t\in\mathbb{R}$. Since
  \begin{linenomath*}
    \begin{equation*}
      \begin{split}
        \int^t_0\tau^5e^{\gamma \tau^2/2}d\tau
        &=\frac{1}{\gamma}t^4\Big(e^{\gamma t^2/2}-1\Big)-\frac{4}{\gamma}\int^t_0\tau^3\Big(e^{\gamma \tau^2/2}-1\Big)d\tau\\
        &\leq \frac{1}{\gamma}t^4\Big(e^{\gamma t^2/2}-1\Big)\qquad\text{for any}~t\geq0,
      \end{split}
    \end{equation*}
  \end{linenomath*}
it follows that
  \begin{linenomath*}
    \begin{equation*}
      |F(t)|\leq \varepsilon t^4+\frac{1}{\gamma}C'_\varepsilon t^4\Big(e^{\gamma t^2/2}-1\Big)\qquad\text{for all}~t\in\mathbb{R}.
    \end{equation*}
  \end{linenomath*}
Also, by the Moser-Trudinger inequality, there exists $C>0$ such that
  \begin{linenomath*}
    \begin{equation*}
      \int_{\mathbb{R}^2}\Big(e^{\gamma u^2}-1\Big)dx\leq C^2\qquad\text{for any}~u\in B_m~\text{with}~\|\nabla u\|_{L^2(\mathbb{R}^2)}\leq 1.
    \end{equation*}
  \end{linenomath*}
Let $\delta\in(0,1)$ be arbitrary. For any $u\in B_m$ with $\|\nabla u\|_{L^2(\mathbb{R}^2)}\leq \delta$,  using also H\"{o}lder inequality and Gagliardo-Nirenberg inequality, we have
  \begin{linenomath*}
    \begin{equation*}
      \begin{split}
        \int_{\mathbb{R}^2}|F(u)|dx
        &\leq \varepsilon \int_{\mathbb{R}^2}u^4dx+\frac{1}{\gamma}C'_\varepsilon\int_{\mathbb{R}^2}u^4\Big(e^{\gamma u^2/2}-1\Big)dx\\
        &\leq \varepsilon \int_{\mathbb{R}^2}u^4dx+\frac{1}{\gamma}C'_\varepsilon\Big(\int_{\mathbb{R}^2}u^8dx\Big)^{\frac{1}{2}}\Big[\int_{\mathbb{R}^2}\Big(e^{\gamma u^2/2}-1\Big)^2dx\Big]^\frac{1}{2}\\
        &\leq \varepsilon \int_{\mathbb{R}^2}u^4dx+\frac{1}{\gamma}C'_\varepsilon\Big(\int_{\mathbb{R}^2}u^8dx\Big)^{\frac{1}{2}}\Big[\int_{\mathbb{R}^2}\left(e^{\gamma u^2}-1\right)dx\Big]^\frac{1}{2}\\
        &\leq \varepsilon C_1m\int_{\mathbb{R}^2}|\nabla u|^2dx+CC_2C'_\varepsilon m^{\frac{1}{2}}(m+1)\Big(\int_{\mathbb{R}^2}|\nabla u|^2dx\Big)^{\frac{3}{2}}\\
        &\leq \Big[\varepsilon C_1 m+CC_2C'_\varepsilon m^{\frac{1}{2}}(m+1)\delta\Big]\int_{\mathbb{R}^2}|\nabla u|^2dx,
      \end{split}
    \end{equation*}
  \end{linenomath*}
where $C_1,C_2>0$ are independent of $m, \varepsilon, \delta$ and $u$. Choosing $\varepsilon>0$ and $\delta\in (0,1)$ small enough, we deduce that \eqref{eq:F1} holds when $N=2$.

Finally we consider the case when $N=1$. Since $H^1(\mathbb{R})\hookrightarrow L^\infty(\mathbb{R})$, there exists $K>0$ such that
  \begin{linenomath*}
    \begin{equation*}
      \|u\|_{L^\infty(\mathbb{R})}\leq K\qquad\text{for any}~u\in B_m~\text{with}~\|\nabla u\|_{L^2(\mathbb{R})}\leq 1.
    \end{equation*}
  \end{linenomath*}
Let $\varepsilon>0$ and $\delta\in (0,1)$ be arbitrary.  By $(f0)$ and $(f1)$, one can find $C''_\varepsilon>0$ such that $|F(t)|\leq \varepsilon t^6+C''_\varepsilon t^{10}$ for all $|t|\leq K$. Thus, for any $u\in B_m$ with $\|\nabla u\|_{L^2(\mathbb{R})}\leq \delta$, it follows that
  \begin{linenomath*}
    \begin{equation*}
      \begin{split}
        \int_{\mathbb{R}}|F(u)|dx
        &\leq \varepsilon \int_{\mathbb{R}}u^6dx+C''_\varepsilon\int_{\mathbb{R}}u^{10}dx\\
        &\leq \varepsilon C_3m^2\int_{\mathbb{R}}|\nabla u|^2dx+C_4C''_\varepsilon m^3\Big(\int_{\mathbb{R}}|\nabla u|^2dx\Big)^2\\
        &\leq\big(\varepsilon C_3 m^2+C_4 C''_\varepsilon m^3 \delta^2\big)\int_{\mathbb{R}}|\nabla u|^2dx,
      \end{split}
    \end{equation*}
  \end{linenomath*}
where $C_3,C_4>0$ are independent of $m, \varepsilon, \delta$ and $u$. Taking $\varepsilon>0$ and $\delta\in(0,1)$ sufficiently small, we derive \eqref{eq:F1} when $N=1$.

$(ii)$ The proofs of the two claims being similar, we only prove that
  \begin{linenomath*}
    \begin{equation}\label{eq:F2}
      \lim_{n\to\infty}\int_{\mathbb{R}^N}\widetilde{F}(u_n)dx=0\qquad\text{if}~\lim_{n\to\infty}\|u_n\|_{L^{2+4/N}(\mathbb{R}^N)}=0.
    \end{equation}
  \end{linenomath*}
When $N=2$, we choose $L>0$ large enough such that $\sup_{n\geq1}\|u_n\|_{H^1(\mathbb{R}^2)}\leq L$. For $\gamma:=1/L^2$, by the Moser-Trudinger inequality, one can find $D>0$ such that
  \begin{linenomath*}
    \begin{equation}\label{eq:MT}
      \sup_{n\geq1}\int_{\mathbb{R}^2}\Big(e^{\gamma u^2_n}-1\Big)dx\leq D.
    \end{equation}
  \end{linenomath*}
Let $\varepsilon>0$ be arbitrary. By $(f0)-(f2)$, there exists $D_\varepsilon>0$ such that $|\widetilde{F}(t)|\leq \varepsilon \big(e^{\gamma t^2}-1\big)+D_\varepsilon t^4$ for all $t\in\mathbb{R}$. Thus
  \begin{linenomath*}
    \begin{equation*}
      \int_{\mathbb{R}^2}|\widetilde{F}(u_n)|dx
      \leq \varepsilon \int_{\mathbb{R}^2}\Big(e^{\gamma u^2_n}-1\Big)dx+D_\varepsilon\int_{\mathbb{R}^2}u^4_ndx\leq \varepsilon D+D_\varepsilon\int_{\mathbb{R}^2}u^4_ndx.
    \end{equation*}
  \end{linenomath*}
Since $\varepsilon$ is arbitrary, it follows that \eqref{eq:F2} holds when $N=2$. The treatment of the cases $N\geq3$ and $N=1$ is similar.

$(iii)$  We only prove \eqref{eq:fuv} when $N=2$ and the other cases follow analogously. Clearly, by \eqref{eq:MT}, we have
  \begin{linenomath*}
    \begin{equation*}
      \sup_{n\geq1}\int_{\mathbb{R}^2}\Big(e^{\gamma u^2_n/2}-1\Big)^2dx\leq \sup_{n\geq1}\int_{\mathbb{R}^2}\Big(e^{\gamma u^2_n}-1\Big)dx\leq D.
    \end{equation*}
  \end{linenomath*}
For given $\varepsilon>0$, the conditions $(f0)-(f2)$ assure the existence of $D'_\varepsilon>0$ such that $|f(t)|\leq \varepsilon \big(e^{\gamma t^2/2}-1\big)+D'_\varepsilon |t|^3$ for all $t\in\mathbb{R}$. Thus
  \begin{linenomath*}
    \begin{equation*}
      \begin{split}
        \int_{\mathbb{R}^2}|f(u_n)v_n|dx
        &\leq \varepsilon \int_{\mathbb{R}^2}\Big(e^{\gamma u^2_n/2}-1\Big)|v_n|dx+D'_\varepsilon\int_{\mathbb{R}^2}|u_n|^3|v_n|dx\\
        &\leq \varepsilon \Big[\int_{\mathbb{R}^2}\left(e^{\gamma u^2_n/2}-1\right)^2dx\Big]^{\frac{1}{2}}\|v_n\|_{L^2(\mathbb{R}^2)}+D'_\varepsilon\|u_n\|^3_{L^4(\mathbb{R}^2)}\|v_n\|_{L^4(\mathbb{R}^2)}\\
        &\leq \varepsilon \sqrt{D}\|v_n\|_{L^2(\mathbb{R}^2)}+D'_\varepsilon\|u_n\|^3_{L^4(\mathbb{R}^2)}\|v_n\|_{L^4(\mathbb{R}^2)}
      \end{split}
    \end{equation*}
  \end{linenomath*}
and we deduce that \eqref{eq:fuv} holds when $N=2$. \hfill $\square$
\begin{remark}\label{remark:P}
Still under the assumptions of Lemma \ref{lemma:I_1}, for any $m>0$, modifying slightly the proof of \eqref{eq:F1}, one can find $\delta=\delta(N,m)>0$ small enough such that
  \begin{linenomath*}
    \begin{equation*}
      \int_{\mathbb{R}^N}|\widetilde{F}(u)|dx\leq \frac{1}{N}\int_{\mathbb{R}^N}|\nabla u|^2dx
    \end{equation*}
  \end{linenomath*}
for all $u\in B_m$ with $\|\nabla u\|_{L^2(\mathbb{R}^N)}\leq \delta$. As a direct consequence,
  \begin{linenomath*}
    \begin{equation*}
      P(u):=\int_{\mathbb{R}^N}|\nabla u|^2dx-\frac{N}{2}\int_{\mathbb{R}^N}\widetilde{F}(u)dx\geq\frac{1}{2}\int_{\mathbb{R}^N}|\nabla u|^2dx
    \end{equation*}
  \end{linenomath*}
for any $u\in B_m$ with $\|\nabla u\|_{L^2(\mathbb{R}^N)}\leq \delta$.
\end{remark}

For any $u\in H^1(\mathbb{R}^N)$ and $s\in\mathbb{R}$, we define the function
  \begin{linenomath*}
    \begin{equation*}
      (s\star u)(x):=e^{Ns/2}u(e^sx) \qquad \text{for a.e.}~x\in\mathbb{R}^N.
    \end{equation*}
  \end{linenomath*}
Clearly, $s\star u\in H^1(\mathbb{R}^N)$ and $\|s\star u\|_{L^2(\mathbb{R}^N)}=\|u\|_{L^2(\mathbb{R}^N)}$ for all $s\in\mathbb{R}$. We fix $u\neq0$ and consider the real valued function $s\mapsto I(s\star u)$ under the conditions $(f0)-(f3)$.
\begin{lemma}\label{lemma:I_2}
  Assume that $N\geq1$ and $f$ satisfies $(f0)-(f3)$. For any $u\in H^1(\mathbb{R}^N)\setminus\{0\}$, one has
      \begin{itemize}
        \item[$(i)$] $I(s\star u)\to 0^+$ as $s\to-\infty$,
        \item[$(ii)$] $I(s\star u)\to-\infty$ as $s\to+\infty$.
      \end{itemize}
\end{lemma}
\proof $(i)$ Let $m:=\|u\|^2_{L^2(\mathbb{R}^N)}>0$. Since $s\star u\in S_m\subset B_m$ and
  \begin{linenomath*}
    \begin{equation*}
      \|\nabla (s\star u)\|_{L^2(\mathbb{R}^N)}=e^s\|\nabla u\|_{L^2(\mathbb{R}^N)},
    \end{equation*}
  \end{linenomath*}
by Lemma \ref{lemma:I_1} $(i)$, it follows that
  \begin{linenomath*}
    \begin{equation*}
      \frac{1}{4}e^{2s}\int_{\mathbb{R}^N}|\nabla u|^2dx\leq I(s\star u)\leq e^{2s}\int_{\mathbb{R}^N}|\nabla u|^2dx\qquad{when}~s\to-\infty.
    \end{equation*}
  \end{linenomath*}
Thus, $\lim_{s\to-\infty}I(s\star u)=0^+$.

$(ii)$ For any $\lambda\geq0$, we define a function $h_\lambda:\mathbb{R}\to\mathbb{R}$ as follows:
  \begin{linenomath*}
    \begin{equation}\label{eq:hlambda}
      h_\lambda(t):=\left\{
             \begin{aligned}
               \frac{F(t)}{|t|^{2+4/N}}+\lambda,&\qquad\text{for}~t\neq0,\\
               \lambda,\qquad&\qquad\text{for}~t=0.
             \end{aligned}
      \right.
    \end{equation}
  \end{linenomath*}
Clearly, $F(t)=h_\lambda(t)|t|^{2+4/N}-\lambda |t|^{2+4/N}$ for all $t\in\mathbb{R}$. Also, from $(f0)$, $(f1)$ and $(f3)$, it follows that $h_\lambda$ is continuous and
  \begin{linenomath*}
    \begin{equation*}
      h_\lambda(t)\to+\infty\qquad\text{as}~t\to\infty.
    \end{equation*}
  \end{linenomath*}
Choose $\lambda>0$ large enough such that $h_\lambda(t)\geq0$ for any $t\in\mathbb{R}$. By Fatou's lemma, we then have
  \begin{linenomath*}
    \begin{equation*}
      \lim_{s\to+\infty}\int_{\mathbb{R}^N}h_\lambda(e^{Ns/2}u)|u|^{2+\frac{4}{N}}dx=+\infty.
    \end{equation*}
  \end{linenomath*}
Since
  \begin{linenomath*}
    \begin{equation}\label{eq:Isu}
      \begin{split}
        I(s\star u)&=\frac{1}{2}\int_{\mathbb{R}^N}|\nabla (s\star u)|^2dx+\lambda\int_{\mathbb{R}^N}|s\star u|^{2+\frac{4}{N}}dx-\int_{\mathbb{R}^N}h_\lambda(s\star u)|s\star u|^{2+\frac{4}{N}}dx\\
        &=e^{2s}\left[\frac{1}{2}\int_{\mathbb{R}^N}|\nabla u|^2dx+\lambda\int_{\mathbb{R}^N}|u|^{2+\frac{4}{N}}dx-\int_{\mathbb{R}^N}h_\lambda(e^{Ns/2}u)|u|^{2+\frac{4}{N}}dx\right],
      \end{split}
    \end{equation}
  \end{linenomath*}
we deduce that $I(s\star u)\to-\infty$ as $s\to+\infty$. \hfill $\square$

We now assume in addition the monotonicity condition $(f4)$ and work out more properties. First we observe
\begin{remark}\label{remark:g}
  Assume that $N\geq1$. If $f$ satisfies $(f0)$, $(f1)$ and $(f4)$, then one can define a continuous function $g:\mathbb{R}\to\mathbb{R}$ as follows:
    \begin{linenomath*}
      \begin{equation}\label{eq:g}
        g(t):=
            \left\{
              \begin{aligned}
                \frac{f(t)t-2F(t)}{|t|^{2+4/N}},&\qquad\text{for}~t\neq0,\\
                0,\qquad\quad&\qquad\text{for}~t=0.
              \end{aligned}
            \right.
      \end{equation}
    \end{linenomath*}
  Moreover, $g$ is strictly decreasing on $(-\infty,0]$ and strictly increasing on $[0,\infty)$.
\end{remark}

\begin{lemma}\label{lemma:fF}
  Assume that $N\geq1$. If $f$ satisfies $(f0)$, $(f1)$, $(f3)$ and $(f4)$,  then
    \begin{linenomath*}
      \begin{equation*}
        f(t)t>\Big(2+\frac{4}{N}\Big)F(t)>0\qquad\text{for all}~t\neq0.
      \end{equation*}
    \end{linenomath*}
\end{lemma}
\proof We split the proof into several claims.

\smallskip

\textbf{Claim 1.} \emph{$F(t)>0$ for any $t\neq0$.}

Indeed, if $F(t_0)\leq0$ for some $t_0\neq0$, by $(f1)$ and $(f3)$, the function $F(t)/|t|^{2+4/N}$ reaches the global minimum at some $\tau\neq0$ satisfying $F(\tau)\leq0$ and
  \begin{linenomath*}
    \begin{equation*}
      \Big[F(t)/|t|^{2+4/N}\Big]'_{t=\tau}=\frac{f(\tau)\tau-\bigl(2+4/N\bigr)F(\tau)}{|\tau|^{3+4/N}\text{sign}(\tau)}=0.
    \end{equation*}
  \end{linenomath*}
Noting that $f(t)t>2F(t)$ for any $t\neq0$ by Remark \ref{remark:g}, we derive a contradiction:
  \begin{linenomath*}
    \begin{equation*}
      0<f(\tau)\tau-2F(\tau)=\frac{4}{N}F(\tau)\leq0.
    \end{equation*}
  \end{linenomath*}
The proof of Claim 1 is complete.

\smallskip

\textbf{Claim 2.} \emph{There exists a positive sequence $\{\tau^+_n\}$ and a negative sequence $\{\tau^-_n\}$ such that $|\tau^\pm_n|\to0$ and $f(\tau^\pm_n)\tau^\pm_n>(2+4/N)F(\tau^\pm_n)$ for each $n\geq1$.}

We first consider the positive case. By contradiction, we assume that there exists $T_s>0$ small enough such that $f(t)t\leq (2+4/N)F(t)$ for any $t\in(0,T_s]$. Using Claim 1, we have
  \begin{linenomath*}
    \begin{equation*}
      F(t)/t^{2+4/N}\geq F(T_s)/T_s^{2+4/N}>0\qquad\text{for all}~t\in(0,T_s].
    \end{equation*}
  \end{linenomath*}
Noting that $\lim_{t\to0}F(t)/|t|^{2+4/N}=0$ by $(f1)$, we obtain a contradiction. The negative case is similar and so we obtain Claim 2.

\smallskip

\textbf{Claim 3.} \emph{There exists a positive sequence $\{\sigma^+_n\}$ and a negative sequence $\{\sigma^-_n\}$ such that $|\sigma^\pm_n|\to +\infty$ and $f(\sigma^\pm_n)\sigma^\pm_n>(2+4/N)F(\sigma^\pm_n)$ for each $n\geq1$.}

The two cases being similar, we only show the existence of $\{\sigma^-_n\}$. Assume by contradiction that there exists $T_l>0$  such that $f(t)t\leq (2+4/N)F(t)$ for any $t\leq-T_l$. We then have
  \begin{linenomath*}
    \begin{equation*}
      F(t)/|t|^{2+4/N}\leq F(-T_l)/T^{2+4/N}_l<+\infty\qquad\text{for all}~t<-T_l,
    \end{equation*}
  \end{linenomath*}
which contradicts $(f3)$. Therefore, the sequence $\{\sigma^-_n\}$ exists and this proves Claim 3.

\smallskip

\textbf{Claim 4.} \emph{$f(t)t\geq(2+4/N)F(t)$ for any $t\neq0$.}

Let us assume by contradiction that $f(t_0)t_0<(2+4/N)F(t_0)$ for some $t_0\neq0$. Since the cases $t_0<0$ and $t_0>0$ can be treated in a similar way, we can assume further that $t_0<0$.  By Claims 2 and 3, there exist $\tau_{min},\tau_{max}\in\mathbb{R}$ such that $\tau_{min}<t_0<\tau_{max}<0$,
  \begin{linenomath*}
    \begin{equation}\label{eq:fF_1}
      f(t)t<(2+4/N)F(t)\qquad\text{for any}~t\in (\tau_{min},\tau_{max}),
    \end{equation}
  \end{linenomath*}
and
  \begin{linenomath*}
    \begin{equation}\label{eq:fF_2}
      f(t)t=(2+4/N)F(t)\qquad\text{when}~t=\tau_{min},\tau_{max}.
    \end{equation}
  \end{linenomath*}
In view of \eqref{eq:fF_1}, we have
  \begin{linenomath*}
    \begin{equation}\label{eq:fF_3}
      \frac{F(\tau_{min})}{|\tau_{\min}|^{2+4/N}}<\frac{F(\tau_{max})}{|\tau_{\max}|^{2+4/N}}.
    \end{equation}
  \end{linenomath*}
On the other hand, by \eqref{eq:fF_2} and $(f4)$, it is clear that
  \begin{linenomath*}
    \begin{equation}\label{eq:fF_4}
      \frac{F(\tau_{min})}{|\tau_{\min}|^{2+4/N}}=\frac{N}{4}\frac{\widetilde{F}(\tau_{min})}{|\tau_{min}|^{2+4/N}}> \frac{N}{4}\frac{\widetilde{F}(\tau_{max})}{|\tau_{max}|^{2+4/N}}=\frac{F(\tau_{max})}{|\tau_{\max}|^{2+4/N}}.
    \end{equation}
  \end{linenomath*}
Since \eqref{eq:fF_3} and \eqref{eq:fF_4} contradict each other, we obtain Claim 4.

\smallskip

\textbf{Claim 5.} \emph{$f(t)t>(2+4/N)F(t)$ for any $t\neq0$.}

By Claim 4, the function $F(t)/|t|^{2+4/N}$ is nonincreasing on $(-\infty,0)$ and nondecreasing on $(0,\infty)$. Then, in view of $(f4)$, the function $f(t)/|t|^{1+4/N}$ is strictly increasing on $(-\infty,0)$ and $(0,\infty)$. For any $t\neq0$, it is clear that
  \begin{linenomath*}
    \begin{equation*}
      \begin{split}
        (2+4/N)F(t)
        &=(2+4/N)\int^t_0f(s)ds\\
        &<(2+4/N)\frac{f(t)}{|t|^{1+4/N}}\int^t_0|s|^{1+4/N}ds=f(t)t
      \end{split}
    \end{equation*}
  \end{linenomath*}
and this proves Claim 5. Now, by Claims 1 and 5, we complete the proof of Lemma \ref{lemma:fF}. \hfill $\square$

Now recall the Pohozaev functional
  \begin{linenomath*}
    \begin{equation*}
      P(u):=\int_{\mathbb{R}^N}|\nabla u|^2dx-\frac{N}{2}\int_{\mathbb{R}^N}\widetilde{F}(u)dx,
    \end{equation*}
  \end{linenomath*}
where $\widetilde{F}(t):=f(t)t-2F(t)$ for any $t\in\mathbb{R}$. As an essential technical result where the monotonicity condition $(f4)$ plays its due role, we have
\begin{lemma}\label{lemma:I_3}
  Assume that $N\geq1$ and $f$ satisfies $(f0)-(f4)$. For any $u\in H^1(\mathbb{R}^N)\setminus\{0\}$, the following statements hold.
    \begin{itemize}
      \item[$(i)$] There exists a unique number $s(u)\in\mathbb{R}$ such that $P(s(u)\star u)=0$.
      \item[$(ii)$] $I(s(u)\star u)>I(s\star u)$ for any $s\neq s(u)$. In particular, $I(s(u)\star u)>0$.
      \item[$(iii)$] The mapping $u\mapsto s(u)$ is continuous in $u\in H^1(\mathbb{R}^N)\setminus\{0\}$.
      \item[$(iv)$] $s(u(\cdot+y))=s(u)$ for any $y\in\mathbb{R}^N$. If $f$ is odd, then one also has $s(-u)=s(u)$.
    \end{itemize}
\end{lemma}
\proof  $(i)$  Since
  \begin{linenomath*}
    \begin{equation*}
      I(s\star u)=\frac{1}{2}e^{2s}\int_{\mathbb{R}^N}|\nabla u|^2dx-e^{-Ns}\int_{\mathbb{R}^N}F(e^{Ns/2}u)dx,
    \end{equation*}
  \end{linenomath*}
we see that $I(s\star u)$ is of class $C^1$ and
  \begin{linenomath*}
    \begin{equation*}
      \frac{d}{ds}I(s\star u)=e^{2s}\int_{\mathbb{R}^N}|\nabla u|^2dx-\frac{N}{2}e^{-Ns}\int_{\mathbb{R}^N}\widetilde{F}(e^{Ns/2}u)dx=P(s\star u).
    \end{equation*}
  \end{linenomath*}
By Lemma \ref{lemma:I_2}, we also have
  \begin{linenomath*}
    \begin{equation*}
     \lim_{s\to-\infty}I(s\star u)=0^+\qquad\text{and}\qquad\lim_{s\to+\infty}I(s\star u)=-\infty.
    \end{equation*}
  \end{linenomath*}
Therefore, $I(s\star u)$ reaches the global maximum at some $s(u)\in\mathbb{R}$ and then
  \begin{linenomath*}
    \begin{equation*}
      P(s(u)\star u)=\frac{d}{ds}I(s(u)\star u)=0.
    \end{equation*}
  \end{linenomath*}
To show the uniqueness, we recall the function $g$ defined by \eqref{eq:g}. Since
  \begin{linenomath*}
    \begin{equation*}
      \widetilde{F}(t)=g(t)|t|^{2+\frac{4}{N}}\qquad\text{for all}~t\in\mathbb{R},
    \end{equation*}
  \end{linenomath*}
it follows that
  \begin{linenomath*}
    \begin{equation*}
      P(s\star u)=e^{2s}\left[\int_{\mathbb{R}^N}|\nabla u|^2dx-\frac{N}{2}\int_{\mathbb{R}^N}g(e^{Ns/2}u)|u|^{2+\frac{4}{N}}dx\right].
    \end{equation*}
  \end{linenomath*}
Noting that, for fixed $t\in\mathbb{R}\setminus\{0\}$, the function $s\mapsto g(e^{Ns/2}t)$ is strictly increasing by $(f4)$ and Remark \ref{remark:g}, we conclude that $s(u)$ is unique.

$(ii)$ This item is a direct consequence of the proof above.

$(iii)$ By Item $(i)$, the mapping $u\mapsto s(u)$ is well-defined. Let $u \in  H^1(\mathbb{R}^N)\setminus\{0\}$ and $\{u_n\}\subset H^1(\mathbb{R}^N)\setminus\{0\}$ be any sequence such that $u_n\to u$ in $H^1(\mathbb{R}^N)$. Setting $s_n:=s(u_n)$ for any $n\geq1$, we only need to prove that up to a subsequence $s_n\to s(u)$ as $n\to\infty$.

We first show that $\{s_n\}$ is bounded. Recall the continuous coercive function $h_\lambda$ defined by \eqref{eq:hlambda}. Clearly, $h_0(t)\geq 0$ for any $t\in\mathbb{R}$ by Lemma \ref{lemma:fF}. If up to a subsequence  $s_n\to+\infty$, by Fatou's lemma and the fact that $u_n\to u\neq 0$ almost everywhere in $\mathbb{R}^N$, we have
  \begin{linenomath*}
    \begin{equation*}
      \lim_{n\to\infty}\int_{\mathbb{R}^N}h_0(e^{Ns_n/2}u_n)|u_n|^{2+\frac{4}{N}}dx=+\infty.
    \end{equation*}
  \end{linenomath*}
In view of Item $(ii)$ and \eqref{eq:Isu} with $\lambda=0$, we then obtain
  \begin{linenomath*}
    \begin{equation}\label{eq:contradiction}
        0\leq e^{-2s_n}I(s_n\star u_n)=\frac{1}{2}\int_{\mathbb{R}^N}|\nabla u_n|^2dx-\int_{\mathbb{R}^N}h_0(e^{Ns_n/2}u_n)|u_n|^{2+\frac{4}{N}}dx\to-\infty,
    \end{equation}
  \end{linenomath*}
which is a contradiction. Therefore, the sequence $\{s_n\}$ is bounded from above. On the other hand, by Item $(ii)$, one has
  \begin{linenomath*}
    \begin{equation*}
      I(s_n\star u_n) \geq I(s(u)\star u_n)\qquad\text{for any}~n\geq1.
    \end{equation*}
  \end{linenomath*}
Since $s(u)\star u_n\to s(u)\star u$ in $H^1(\mathbb{R}^N)$, it follows that
  \begin{linenomath*}
    \begin{equation*}
      I(s(u)\star u_n)= I(s(u)\star u)+o_n(1)
    \end{equation*}
  \end{linenomath*}
and thus
  \begin{linenomath*}
    \begin{equation}\label{eq:liminf}
      \liminf_{n\to\infty}I(s_n\star u_n)\geq I(s(u)\star u)>0.
    \end{equation}
  \end{linenomath*}
As $\{s_n \star u_n\}\subset B_m$ for $m>0$ large enough, in view of Lemma \ref{lemma:I_1} $(i)$ and the fact that
  \begin{linenomath*}
    \begin{equation*}
      \|\nabla (s_n\star u_n)\|_{L^2(\mathbb{R}^N)}=e^{s_n}\|\nabla u_n\|_{L^2(\mathbb{R}^N)},
    \end{equation*}
  \end{linenomath*}
we deduce from \eqref{eq:liminf} that $\{s_n\}$ is bounded also from below.

Without loss of generality, we can now assume that
  \begin{linenomath*}
    \begin{equation*}
      s_n\to s_*\qquad\text{for some}~s_*\in\mathbb{R}.
    \end{equation*}
  \end{linenomath*}
Recalling that $u_n\to u$ in $H^1(\mathbb{R}^N)$, one then has $s_n\star u_n\to s_*\star u$ in $H^1(\mathbb{R}^N)$. Since $P(s_n\star u_n)=0$ for any $n\geq1$, it follows that
  \begin{linenomath*}
    \begin{equation*}
      P(s_*\star u)=0.
    \end{equation*}
  \end{linenomath*}
By Item $(i)$, we see that $s_*=s(u)$ and thus Item $(iii)$ is proved.

$(iv)$ For any $y\in\mathbb{R}^N$, by changing variables in the integrals, we have
  \begin{linenomath*}
    \begin{equation*}
      P\bigl(s(u)\star u(\cdot+y)\bigr)=P\bigl(s(u)\star u\bigr)=0
    \end{equation*}
  \end{linenomath*}
and thus $s(u(\cdot+y))=s(u)$ via Item $(i)$. When $f$ is odd, it is clear that
  \begin{linenomath*}
    \begin{equation*}
      P\bigl(s(u)\star (-u)\bigr)=P\bigl(-(s(u)\star u)\bigr)=P\bigl(s(u)\star u\bigr)=0
    \end{equation*}
  \end{linenomath*}
and hence $s(-u)=s(u)$. \hfill $\square$

Under the assumptions of Lemma \ref{lemma:I_3}, we also have the following result which concerns the Pohozaev manifold
  \begin{linenomath*}
    \begin{equation*}
      \mathcal{P}_m:=\Big\{u\in S_m~\Big|~P(u)=0\Big\}
    \end{equation*}
  \end{linenomath*}
and the functional $I$ constrained to it.
\begin{lemma}\label{lemma:I_4}
Assume that $N \geq 1$ and $f$ satisfies $(f0)-(f4)$. Then
  \begin{itemize}
    \item[$(i)$] $\mathcal{P}_m\neq\emptyset$,
    \item[$(ii)$] $\inf_{u\in \mathcal{P}_m}\|\nabla u\|_{L^2(\mathbb{R}^N)}>0$,
    \item[$(iii)$] $\inf_{u\in\mathcal{P}_m}I(u)>0$,
    \item[$(iv)$] $I$ is coercive on $\mathcal{P}_m$, that is $I(u_n)\to+\infty$ for any $\{u_n\}\subset\mathcal{P}_m$ with $\|u_n\|_{H^1(\mathbb{R}^N)}\to\infty$.
  \end{itemize}
\end{lemma}
\proof $(i)$ This item is a direct consequence of Lemma \ref{lemma:I_3} $(i)$.

$(ii)$ If there exists $\{u_n\}\subset \mathcal{P}_m$ such that $\|\nabla u_n\|_{L^2(\mathbb{R}^N)}\to0$, then Remark \ref{remark:P} implies
  \begin{linenomath*}
    \begin{equation*}
      0=P(u_n)\geq\frac{1}{2}\int_{\mathbb{R}^N}|\nabla u_n|^2dx>0\qquad\text{for}~n~\text{large enough},
    \end{equation*}
  \end{linenomath*}
which is a contradiction. Therefore, $\inf_{u\in \mathcal{P}_m}\|\nabla u\|_{L^2(\mathbb{R}^N)}>0$.

$(iii)$ For any $u\in\mathcal{P}_m$, by Lemma \ref{lemma:I_3} $(i)$ and $(ii)$, we have
  \begin{linenomath*}
    \begin{equation*}
      I(u)=I(0\star u)\geq I(s\star u)\qquad\text{for all}~s\in \mathbb{R}.
    \end{equation*}
  \end{linenomath*}
Let $\delta>0$ be the number given by Lemma \ref{lemma:I_1} $(i)$ and $s:=\ln \big(\delta/\|\nabla u\|_{L^2(\mathbb{R}^N)}\big)$. Since $\|\nabla (s\star u)\|_{L^2(\mathbb{R}^N)}=\delta$, by Lemma \ref{lemma:I_1} $(i)$, we deduce that
  \begin{linenomath*}
    \begin{equation*}
      I(u)\geq I(s\star u)\geq\frac{1}{4}\int_{\mathbb{R}^N}|\nabla (s\star u)|^2dx=\frac{1}{4}\delta^2
    \end{equation*}
  \end{linenomath*}
and thus Item $(iii)$ holds.

$(iv)$ By contradiction, we assume that there exists $\{u_n\}\subset \mathcal{P}_m$ such that $\|u_n\|_{H^1(\mathbb{R}^N)}\to\infty$ but $\sup_{n\geq1}I(u_n)\leq c$ for some $c\in(0,+\infty)$. For any $n\geq1$, set
  \begin{linenomath*}
    \begin{equation*}
      s_n:=\ln\bigl(\|\nabla u_n\|_{L^2(\mathbb{R}^N)}\bigr)\qquad\text{and}\qquad v_n:=(-s_n)\star u_n.
    \end{equation*}
  \end{linenomath*}
Clearly, $s_n\to+\infty$, $\{v_n\}\subset S_m$ and $\|\nabla v_n\|_{L^2(\mathbb{R}^N)}=1$ for any $n\geq1$. Let
  \begin{linenomath*}
    \begin{equation*}
      \rho:=\limsup_{n\to\infty}\Big(\sup_{y\in\mathbb{R}^N}\int_{B(y,1)}|v_n|^2dx\Big).
    \end{equation*}
  \end{linenomath*}
To derive a contradiction, we distinguish the two cases: \emph{non-vanishing} and \emph{vanishing}.

$\bullet$  \emph{Non-vanishing: that is $\rho>0$}. Up to a subsequence, there exists $\{y_n\}\subset \mathbb{R}^N$ and $w\in H^1(\mathbb{R}^N)\setminus\{0\}$ such that
  \begin{linenomath*}
    \begin{equation*}
      w_n:=v_n(\cdot+y_n)\rightharpoonup w\quad\text{in}~H^1(\mathbb{R}^N)\qquad\text{and}\qquad w_n\to w\quad\text{a.e. in} ~\mathbb{R}^N.
    \end{equation*}
  \end{linenomath*}
Recall the continuous coercive function $h_\lambda$ defined by \eqref{eq:hlambda} and let $\lambda=0$. Since $s_n\to+\infty$, by Lemma \ref{lemma:fF} and Fatou's lemma, it follows that
  \begin{linenomath*}
    \begin{equation*}
      \lim_{n\to\infty}\int_{\mathbb{R}^N}h_0(e^{Ns_n/2}w_n)|w_n|^{2+\frac{4}{N}}dx=+\infty.
    \end{equation*}
  \end{linenomath*}
In view of Item $(iii)$ and \eqref{eq:Isu} with $\lambda=0$, we have
  \begin{linenomath*}
    \begin{equation*}
      \begin{split}
        0 \leq e^{-2s_n}I(u_n)
        &=e^{-2s_n}I(s_n\star v_n)\\
        &=\frac{1}{2}-\int_{\mathbb{R}^N}h_0(e^{Ns_n/2}v_n)|v_n|^{2+\frac{4}{N}}dx\\
        &=\frac{1}{2}-\int_{\mathbb{R}^N}h_0(e^{Ns_n/2}w_n)|w_n|^{2+\frac{4}{N}}dx\to-\infty,
      \end{split}
    \end{equation*}
  \end{linenomath*}
which is a contradiction.

$\bullet$  \emph{Vanishing: that is $\rho=0$}. In this case, by Lions Lemma \cite[Lemma I.1]{Lions84-2}, $v_n\to 0$ in $L^{2+4/N}(\mathbb{R}^N)$. Using Lemma \ref{lemma:I_1} $(ii)$, we thus have
  \begin{linenomath*}
    \begin{equation*}
      \lim_{n\to\infty}e^{-Ns}\int_{\mathbb{R}^N}F(e^{Ns/2}v_n)dx=0\qquad\text{for any}~s\in\mathbb{R}.
    \end{equation*}
  \end{linenomath*}
Since $P(s_n\star v_n)=P(u_n)=0$, by Lemma \ref{lemma:I_3} $(i)$ and $(ii)$ it follows that, for any $s \in \mathbb{R}$,
  \begin{linenomath*}
    \begin{equation*}
      \begin{split}
        c\geq I(u_n)=I(s_n\star v_n)
        &\geq I(s\star v_n)\\
        &=\frac{1}{2}e^{2s}-e^{-Ns}\int_{\mathbb{R}^N}F(e^{Ns/2}v_n)dx=\frac{1}{2}e^{2s}+o_n(1).
      \end{split}
    \end{equation*}
  \end{linenomath*}
Clearly, this leads a contradiction for $s>\ln(2c)/2$. Therefore, $I$ is coercive on $\mathcal{P}_m$. \hfill $\square$
\begin{remark}\label{remark:boundedness}
  Assume that $N\geq1$ and $f$ satisfies $(f0)-(f4)$. For any sequence $\{u_n\}\subset H^1(\mathbb{R}^N)\setminus\{0\}$ such that
    \begin{linenomath*}
      \begin{equation*}
        P(u_n)=0,\qquad \sup_{n\geq1}\|u_n\|_{L^2(\mathbb{R}^N)}<+\infty\quad\text{and}\quad \sup_{n\geq1}I(u_n)<+\infty.
      \end{equation*}
    \end{linenomath*}
  repeating the proof of Lemma \ref{lemma:I_4} $(iv)$, one has that $\{u_n\}$ is bounded in $H^1(\mathbb{R}^N)$.
\end{remark}

To end this section, we give a Brezis-Lieb type splitting result which is needed when we study the convergence of the Palais-Smale sequences.
\begin{lemma}\label{lemma:BL}
  Assume that $N\geq1$ and  $f$ is a continuous function satisfying the conditions below:
    \begin{itemize}
      \item[$(C1)$] when $N=1$, for any $T>0$, there exists $C_T>0$ such that $|f(t)|\leq C_T|t|$ for all $|t|\leq T$;
      \item[$(C2)$] when $N=2$, for any $\gamma>0$, there exists $C_\gamma>0$ such that
                     \begin{linenomath*}
                       \begin{equation*}
                         |f(t)|\leq C_\gamma \Big[|t|+\Big(e^{\gamma t^2}-1\Big)\Big]\qquad\text{for all}~t\in\mathbb{R};
                       \end{equation*}
                     \end{linenomath*}
      \item[$(C3)$] when $N\geq3$, there exists $C>0$ such that $|f(t)|\leq C\left(|t|+|t|^{2^*-1}\right)$ for all $t\in\mathbb{R}$.
    \end{itemize}
  If $\{u_n\}\subset H^1(\mathbb{R}^N)$ is bounded and $u_n\to u$ almost everywhere in $\mathbb{R}^N$ for some $u\in H^1(\mathbb{R}^N)$, then
    \begin{linenomath*}
      \begin{equation}\label{eq:BL}
        \lim_{n\to\infty}\int_{\mathbb{R}^N}\big|F(u_n)-F(u_n-u)-F(u)\big|dx=0.
      \end{equation}
    \end{linenomath*}
\end{lemma}
\proof For the case $N\geq2$, one can find a detailed proof in \cite[Lemma 3.2]{JL19}. Here we only prove \eqref{eq:BL} when $N=1$. Since $H^1(\mathbb{R})\hookrightarrow L^\infty(\mathbb{R})$, there exists $T>0$ large enough such that
  \begin{linenomath*}
    \begin{equation*}
      \sup_{n\geq1}\|u_n\|_{L^\infty(\mathbb{R})},~\sup_{n\geq1}\|u_n-u\|_{L^\infty(\mathbb{R})},~\|u\|_{L^\infty(\mathbb{R})}\leq T.
    \end{equation*}
  \end{linenomath*}
Let $\varepsilon>0$ be arbitrary. For any $a,b\in\mathbb{R}$ such that $|a|,|b|,|a+b|\leq T$, by $(C1)$ and Young's inequality, we have
  \begin{linenomath*}
    \begin{equation*}
      \begin{split}
        |F(a+b)-F(a)|&=\Big|\int^1_0f(a+\tau b)bd\tau\Big|\\
        &\leq C_T\int^1_0|a+\tau b||b|d\tau\leq C_T(|a||b|+b^2)\\
        &\leq \varepsilon C_T a^2+C_T(1+\varepsilon^{-1})b^2=:\varepsilon \varphi(a)+\psi_\varepsilon(b).
      \end{split}
    \end{equation*}
  \end{linenomath*}
In particular, $|F(b)|\leq \psi_\varepsilon(b)$ for any $|b|\leq T$. Note that $\int_{\mathbb{R}}\varphi(u_n-u)dx$ is bounded uniformly in $\varepsilon$ and $n$, $\int_\mathbb{R}\psi_\varepsilon(u)dx<\infty$ for any $\varepsilon>0$, and $F(u)\in L^1(\mathbb{R})$. Repeating the argument in \cite[Proof of Theorem 2]{BL83}, we deduce that \eqref{eq:BL} holds when $N=1$. \hfill $\square$

\section{The behavior of the function $m \mapsto E_m$}\label{sect:Em}

When $N\geq1$ and $f$ satisfies $(f0)-(f4)$, for given $m>0$, one can see from Lemma \ref{lemma:I_4} that the infimum
  \begin{linenomath*}
    \begin{equation*}
      E_m:=\inf_{u \in \mathcal{P}_m} I(u)
    \end{equation*}
  \end{linenomath*}
is well defined and strictly positive. Our goal in this section is to characterize further the behavior of $E_m$ when $m > 0$ varies. In particular we shall prove that $E_m$ is nonincreasing in $m>0$.  We start by showing the continuity of $E_m$.
\begin{lemma}\label{lemma:Em_continuity}
Assume that $N \geq 1$ and $f$ satisfies $(f0)-(f4)$. Then the function $m\mapsto E_m$ is continuous at each $m>0$.
\end{lemma}
\proof It is equivalent to prove that for a given $m>0$ and any positive sequence $\{m_k\}$ such that $m_k\to m$ as $k\to\infty$, one has $\lim_{k\to\infty}E_{m_k}=E_m$. We first show that
  \begin{linenomath*}
    \begin{equation}\label{eq:Em1_1}
      \limsup_{k\to\infty}E_{m_k}\leq E_m.
    \end{equation}
  \end{linenomath*}
For any $u\in \mathcal{P}_m$, we define
  \begin{linenomath*}
    \begin{equation*}
      u_k:=\sqrt{\frac{m_k}{m}}u\in S_{m_k},\qquad k\in\mathbb{N}^+.
    \end{equation*}
  \end{linenomath*}
Since $u_k\to u$ in $H^1(\mathbb{R}^N)$, by Lemma \ref{lemma:I_3} $(iii)$, we have $\lim_ {k\to\infty}s(u_k)=s(u)=0$ and thus
  \begin{linenomath*}
    \begin{equation*}
      s(u_k)\star u_k\to s(u)\star u=u\quad\text{in}~H^1(\mathbb{R}^N)~\text{as}~k\to\infty.
    \end{equation*}
  \end{linenomath*}
As a consequence,
  \begin{linenomath*}
    \begin{equation*}
      \limsup_{k\to\infty}E_{m_k}\leq \limsup_{k\to\infty} I\bigl(s(u_k)\star u_k\bigr)=I(u).
    \end{equation*}
  \end{linenomath*}
Noting that $u\in\mathcal{P}_m$ is arbitrary, we deduce that \eqref{eq:Em1_1} holds.

To complete the proof, it remains to show that
  \begin{linenomath*}
    \begin{equation}\label{eq:Em1_2}
      \liminf_{k\to\infty}E_{m_k}\geq E_m.
    \end{equation}
  \end{linenomath*}
For each $k\in\mathbb{N}^+$, there exists $v_k\in\mathcal{P}_{m_k}$ such that
  \begin{linenomath*}
    \begin{equation}\label{eq:Em1_3}
      I(v_k)\leq E_{m_k}+\frac{1}{k}.
    \end{equation}
  \end{linenomath*}
Setting
  \begin{linenomath*}
    \begin{equation*}
      t_k:=\Big(\frac{m}{m_k}\Big)^{1/N}\qquad\text{and}\qquad\tilde{v}_k:=v_k(\cdot/t_k)\in S_m,
    \end{equation*}
  \end{linenomath*}
by Lemma \ref{lemma:I_3} $(ii)$ and \eqref{eq:Em1_3}, we have
  \begin{linenomath*}
    \begin{equation*}
      \begin{split}
        E_m\leq I\bigl(s(\tilde{v}_k)\star \tilde{v}_k\bigr)
        &\leq I\bigl(s(\tilde{v}_k)\star v_k\bigr)+\Big|I\bigl(s(\tilde{v}_k)\star \tilde{v}_k\bigr)-I\bigl(s(\tilde{v}_k)\star v_k\bigr)\Big|\\
        &\leq I(v_k)+\Big|I\bigl(s(\tilde{v}_k)\star \tilde{v}_k\bigr)-I\bigl(s(\tilde{v}_k)\star v_k\bigr)\Big|\\
        &\leq E_{m_k}+\frac{1}{k}+\Big|I\bigl(s(\tilde{v}_k)\star \tilde{v}_k\bigr)-I\bigl(s(\tilde{v}_k)\star v_k\bigr)\Big|\\
        &=:E_{m_k}+\frac{1}{k}+C(k).
      \end{split}
    \end{equation*}
  \end{linenomath*}
It is clear that one will obtain \eqref{eq:Em1_2} if
  \begin{linenomath*}
    \begin{equation}\label{eq:Em1_4}
      \lim_{k\to\infty}C(k)=0.
    \end{equation}
  \end{linenomath*}
Noting that $s\star (v(\cdot/t))=(s\star v)(\cdot/t)$, we have
  \begin{linenomath*}
    \begin{equation*}
      \begin{split}
        C(k)
        &=\left|\frac{1}{2}\bigl(t^{N-2}_k-1\bigr)\int_{\mathbb{R}^N}\bigl|\nabla \bigl(s(\tilde{v}_k)\star v_k\bigr) \bigr|^2dx-\bigl(t^N_k-1\bigr)\int_{\mathbb{R}^N}F\bigl(s(\tilde{v}_k)\star v_k\bigr)dx\right|\\
        &\leq \frac{1}{2}\bigl|t^{N-2}_k-1\bigr|\cdot\int_{\mathbb{R}^N}\bigl|\nabla \bigl(s(\tilde{v}_k)\star v_k\bigr) \bigr|^2dx+\bigl|t^N_k-1\bigr|\cdot\int_{\mathbb{R}^N}\bigl|F\bigl(s(\tilde{v}_k)\star v_k\bigr)\bigr|dx\\
        &=:\frac{1}{2}\bigl|t^{N-2}_k-1\bigr|\cdot A(k)+\bigl|t^N_k-1\bigr|\cdot B(k).
      \end{split}
    \end{equation*}
  \end{linenomath*}
Since $t_k\to1$, the proof of \eqref{eq:Em1_4} and thus of \eqref{eq:Em1_2} is reduced to showing that
  \begin{linenomath*}
    \begin{equation}\label{eq:Em1_5}
      \limsup_{k\to\infty}A(k)<+\infty\qquad\text{and}\qquad \limsup_{k\to\infty}B(k)<+\infty.
    \end{equation}
  \end{linenomath*}
To justify \eqref{eq:Em1_5}, we prove below three claims in turn.

\smallskip

\textbf{Claim 1.} \emph{The sequence $\{v_k\}$ is bounded in $H^1(\mathbb{R}^N)$.}

Indeed, by \eqref{eq:Em1_3} and \eqref{eq:Em1_1},
  \begin{linenomath*}
    \begin{equation*}
      \limsup_{k\to\infty}I(v_k)\leq E_m.
    \end{equation*}
  \end{linenomath*}
Since $v_k\in\mathcal{P}_{m_k}$ and $m_k\to m$, we deduce from Remark \ref{remark:boundedness} that Claim 1 holds.

\smallskip

\textbf{Claim 2.} \emph{The sequence $\{\tilde{v}_k\}$ is bounded in $H^1(\mathbb{R}^N)$, and there exists $\{y_k\}\subset\mathbb{R}^N$ and $v\in H^1(\mathbb{R}^N)$ such that up to a subsequence $\tilde{v}_k(\cdot+y_k)\to v\neq0$ almost everywhere in $\mathbb{R}^N$.}

Indeed, since $t_k\to1$, it follows from Claim 1 that $\{\tilde{v}_k\}$ is bounded in $H^1(\mathbb{R}^N)$. Set
  \begin{linenomath*}
    \begin{equation*}
      \rho:=\limsup_{k\to\infty}\Big(\sup_{y\in\mathbb{R}^N}\int_{B(y,1)}|\tilde{v}_k|^2dx\Big).
    \end{equation*}
  \end{linenomath*}
We now only need to rule out the case $\rho=0$. If $\rho=0$, then $\tilde{v}_k\to0$ in $L^{2+4/N}(\mathbb{R}^N)$ by Lions Lemma \cite[Lemma I.1]{Lions84-2}. As a consequence,
  \begin{linenomath*}
    \begin{equation*}
      \int_{\mathbb{R}^N}|v_k|^{2+\frac{4}{N}}dx
      =\int_{\mathbb{R}^N}|\tilde{v}_k(t_k\cdot)|^{2+\frac{4}{N}}dx=t^{-N}_k\int_{\mathbb{R}^N}|\tilde{v}_k|^{2+\frac{4}{N}}dx\to0.
    \end{equation*}
  \end{linenomath*}Combining Lemma \ref{lemma:I_1} $(ii)$ and that $P(v_k)=0$, we have
  \begin{linenomath*}
    \begin{equation*}
      \int_{\mathbb{R}^N}|\nabla v_k|^2dx=\frac{N}{2}\int_{\mathbb{R}^N}\widetilde{F}(v_k)dx\to0.
    \end{equation*}
  \end{linenomath*}
In view of Remark \ref{remark:P}, we thus obtain
  \begin{linenomath*}
    \begin{equation*}
      0=P(v_k)\geq\frac{1}{2}\int_{\mathbb{R}^N}|\nabla v_k|^2dx>0\qquad\text{for}~k~\text{large enough},
    \end{equation*}
  \end{linenomath*}
which is a contradiction. The proof of Claim 2 is complete.

\smallskip

\textbf{Claim 3.} \emph{$\limsup_{k\to\infty} s(\tilde{v}_k)<+\infty$.}

Indeed, if Claim 3 does not hold, then up to a subsequence
  \begin{linenomath*}
    \begin{equation}\label{eq:Em1_6}
      s(\tilde{v}_k)\to+\infty\qquad\text{as}~k\to\infty.
    \end{equation}
  \end{linenomath*}
To derive a contradiction, we make some observations at first. By Claim 2, we see that up to a subsequence
  \begin{linenomath*}
    \begin{equation}\label{eq:Em1_7}
      \tilde{v}_k(\cdot+y_k)\to v\neq0\qquad\text{a.e. in}~\mathbb{R}^N.
    \end{equation}
  \end{linenomath*}
On the other hand, Lemma \ref{lemma:I_3} $(iv)$ and \eqref{eq:Em1_6} imply
  \begin{linenomath*}
    \begin{equation}\label{eq:Em1_8}
      s(\tilde{v}_k(\cdot+y_k))=s(\tilde{v}_k)\to+\infty,
    \end{equation}
  \end{linenomath*}
and Lemma \ref{lemma:I_3} $(ii)$ gives us that
  \begin{linenomath*}
    \begin{equation}\label{eq:Em1_9}
      I\bigl(s(\tilde{v}_k(\cdot+y_k))\star \tilde{v}_k(\cdot+y_k)\bigr)\geq0.
    \end{equation}
  \end{linenomath*}
Now, using \eqref{eq:Em1_7}, \eqref{eq:Em1_8} and \eqref{eq:Em1_9}, we clearly obtain a contradiction in the same way to the derivation of \eqref{eq:contradiction}. The proof of Claim 3 is complete.

Now, by Claims 1 and 3, we have
  \begin{linenomath*}
    \begin{equation*}
      \limsup_{k\to\infty}\bigl\|s(\tilde{v}_k)\star v_k\bigr\|_{H^1(\mathbb{R}^N)}<+\infty.
    \end{equation*}
  \end{linenomath*}
Since $f$ satisfies $(f0)-(f2)$, it is clear that \eqref{eq:Em1_5} holds and the lemma is proved. \hfill $\square$

\begin{lemma}\label{lemma:Em_nonincreasing}
Assume that $N \geq 1$ and $f$ satisfies $(f0)-(f4)$. Then the function $m\mapsto E_m$ is nonincreasing on $(0, \infty)$.
\end{lemma}
\proof We only need to show that for any $m>m'>0$ and any arbitrary $\varepsilon>0$ one has
  \begin{linenomath*}
    \begin{equation}\label{eq:Em2_1}
      E_m\leq E_{m'}+\varepsilon.
    \end{equation}
  \end{linenomath*}
By the definition of $E_{m'}$, there exists $u\in\mathcal{P}_{m'}$ such that
  \begin{linenomath*}
    \begin{equation}\label{eq:Em2_2}
      I(u)\leq E_{m'}+\frac{\varepsilon}{2}.
    \end{equation}
  \end{linenomath*}
Let $\chi\in C^\infty_0(\mathbb{R}^N)$ be radial and such that
  \begin{linenomath*}
    \begin{equation*}
      \chi(x)=
        \left\{
               \begin{aligned}
                 &1, &&|x|\leq1,\\
                 &\in[0,1], &&|x|\in(1,2),\\
                 &0, &&|x|\geq2.
               \end{aligned}
        \right.
    \end{equation*}
  \end{linenomath*}
For any small $\delta>0$, we define $u_\delta(x)=u(x)\cdot\chi(\delta x)\in H^1(\mathbb{R}^N)\setminus\{0\}$. Since $u_\delta\to u$ in $H^1(\mathbb{R}^N)$ as $\delta\to0^+$, by Lemma \ref{lemma:I_3} $(iii)$, one has $\lim_{\delta\to0^+}s(u_\delta)=s(u)=0$ and thus
  \begin{linenomath*}
    \begin{equation*}
      s(u_\delta)\star u_\delta\to s(u)\star u=u\quad\text{in}~H^1(\mathbb{R}^N)~\text{as}~\delta\to0^+.
    \end{equation*}
  \end{linenomath*}
As a consequence, we can fix a  $\delta>0$ small enough such that
  \begin{linenomath*}
    \begin{equation}\label{eq:Em2_3}
      I(s(u_\delta)\star u_\delta)\leq I(u)+\frac{\varepsilon}{4}.
    \end{equation}
  \end{linenomath*}
Now take $v\in C^\infty_0(\mathbb{R}^N)$ such that $\text{supp}(v)\subset B(0,1+4/\delta)\setminus B(0,4/\delta)$ and set
  \begin{linenomath*}
    \begin{equation*}
      \tilde{v}=\frac{m-\|u_\delta\|^2_{L^2(\mathbb{R}^N)}}{\|v\|^2_{L^2(\mathbb{R}^N)}}v.
    \end{equation*}
  \end{linenomath*}
For any $\lambda\leq 0$, we define $w_\lambda=u_\delta+\lambda\star \tilde{v}$. Noting that
  \begin{linenomath*}
    \begin{equation*}
      \text{supp}(u_\delta)\cap \text{supp}(\lambda\star \tilde{v})=\emptyset,
    \end{equation*}
  \end{linenomath*}
one has $w_\lambda\in S_m$. We claim that $s(w_\lambda)$ is bounded from above when $\lambda\to-\infty$. Indeed, observing that $I(s(w_\lambda)\star w_\lambda)\geq 0$ by Lemma \ref{lemma:I_3} $(ii)$ and that $w_{\lambda}\to u_\delta\neq0$ almost everywhere in $\mathbb{R}^N$ as $\lambda \to - \infty$, we obtain a contradiction in the same way to the derivation of \eqref{eq:contradiction} if the claim above does not hold. Now since
  \begin{linenomath*}
    \begin{equation*}
      s(w_\lambda)+\lambda\to-\infty\qquad\text{as}~\lambda\to-\infty,
    \end{equation*}
  \end{linenomath*}
we have
  \begin{linenomath*}
    \begin{equation*}
      \bigl\|\nabla[(s(w_\lambda)+\lambda)\star \tilde{v}]\bigr\|_{L^2(\mathbb{R}^N)}\to0\qquad\text{and}\qquad(s(w_\lambda)+\lambda)\star \tilde{v}\to 0\quad\text{in}~L^{2+4/N}(\mathbb{R}^N).
    \end{equation*}
  \end{linenomath*}
By Lemma \ref{lemma:I_1} $(ii)$, it then follows that
  \begin{linenomath*}
    \begin{equation}\label{eq:Em2_4}
      I\bigl((s(w_\lambda)+\lambda)\star\tilde{v}\bigr)\leq\frac{\varepsilon}{4}\qquad\text{for}~\lambda<0~\text{small enough}.
    \end{equation}
  \end{linenomath*}
Now, using Lemma \ref{lemma:I_3} $(ii)$, \eqref{eq:Em2_4}, \eqref{eq:Em2_3} and \eqref{eq:Em2_2}, we obtain
  \begin{linenomath*}
    \begin{equation*}
      \begin{split}
        E_m\leq I\bigl(s(w_\lambda)\star w_\lambda\bigr)
        &=I\bigl(s(w_\lambda)\star u_\delta\bigr)+I\bigl(s(w_\lambda)\star (\lambda\star \tilde{v})\bigr)\\
        &\leq I(s(u_\delta)\star u_\delta)+I\bigl((s(w_\lambda)+\lambda)\star \tilde{v}\bigr)\\
        &\leq I(u)+\frac{\varepsilon}{2}\leq E_{m'}+\varepsilon,
      \end{split}
    \end{equation*}
  \end{linenomath*}
that is \eqref{eq:Em2_1}. \hfill $\square$

\begin{lemma}\label{lemma:Em_strictlymonotonic}
Assume that
$N \geq 1$ and $f$ satisfies $(f0)-(f4)$. Suppose that there exists $u\in S_m$ and $\mu\in\mathbb{R}$ such that
    \begin{linenomath*}
      \begin{equation*}
        -\Delta u+\mu u=f(u)
      \end{equation*}
    \end{linenomath*}
  and $I(u)=E_m$. Then $E_m>E_{m'}$ for any $m'>m$ close enough to $m$ if $\mu>0$ and  for each $m'<m$ near enough to $m$ if $\mu<0$.
\end{lemma}
\proof For any $t>0$ and $s\in\mathbb{R}$, we set $u_{t,s}:=s\star (tu)\in S_{mt^2}$. Since
  \begin{linenomath*}
    \begin{equation*}
      \alpha(t,s):=I(u_{t,s})=\frac{1}{2}t^2e^{2s}\int_{\mathbb{R}^N}|\nabla u|^2dx-e^{-Ns}\int_{\mathbb{R}^N}F(te^{Ns/2}u)dx,
    \end{equation*}
  \end{linenomath*}
it is clear that
  \begin{linenomath*}
    \begin{equation*}
        \frac{\partial}{\partial t}\alpha(t,s)
        =te^{2s}\int_{\mathbb{R}^N}|\nabla u|^2dx-e^{-Ns}\int_{\mathbb{R}^N}f(te^{Ns/2}u)e^{Ns/2}udx=t^{-1}I'(u_{t,s})u_{t,s}.
    \end{equation*}
  \end{linenomath*}
When $\mu>0$, combining the facts that $u_{t,s}\to u$ in $H^1(\mathbb{R}^N)$ as $(t,s)\to(1,0)$ and that
  \begin{linenomath*}
    \begin{equation*}
      I'(u)u=-\mu \|u\|^2_{L^2(\mathbb{R}^N)}=-\mu m<0,
    \end{equation*}
  \end{linenomath*}
one can fix a $\delta>0$ small enough such that
  \begin{linenomath*}
    \begin{equation*}
      \frac{\partial}{\partial t}\alpha(t,s)<0\qquad\text{for any}~(t,s)\in(1,1+\delta]\times[-\delta,\delta].
    \end{equation*}
  \end{linenomath*}
From the mean value theorem, we then obtain
  \begin{linenomath*}
    \begin{equation}\label{eq:Em3_1}
      \alpha(t,s)=\alpha(1,s)+(t-1)\cdot\frac{\partial}{\partial t}\alpha(\theta,s)<\alpha(1,s),
    \end{equation}
  \end{linenomath*}
where $1<\theta<t\leq 1+\delta$ and $|s|\leq\delta$. Note that $s(tu)\to s(u)=0$ as $t\to1^+$ by Lemma \ref{lemma:I_3} $(iii)$. For any $m'>m$ close enough to $m$, we have
  \begin{linenomath*}
    \begin{equation*}
      t:=\sqrt{\frac{m'}{m}}\in(1,1+\delta]\qquad\text{and}\qquad s:=s(tu)\in[-\delta,\delta]
    \end{equation*}
  \end{linenomath*}
and thus, using \eqref{eq:Em3_1} and Lemma \ref{lemma:I_3} $(ii)$,
  \begin{linenomath*}
    \begin{equation*}
      E_{m'}\leq \alpha(t,s(tu))<\alpha(1,s(tu))=I(s(tu)\star u)\leq I(u)=E_m.
    \end{equation*}
  \end{linenomath*}
The case $\mu<0$ can be proved similarly. \hfill $\square$

At this point, from Lemmas \ref{lemma:Em_nonincreasing} and \ref{lemma:Em_strictlymonotonic} we directly obtain
\begin{lemma}\label{lemma:Em_decreasing}
  Assume that $N\geq1$ and $f$ satisfies $(f0)-(f4)$. If there exists $u\in S_m$ and $\mu\in\mathbb{R}$ such that
    \begin{linenomath*}
      \begin{equation*}
        -\Delta u+\mu u=f(u)
      \end{equation*}
    \end{linenomath*}
  and $I(u)=E_m$, then $\mu\geq0$. If in addition $\mu>0$, then $E_m>E_{m'}$ for any $m'>m$.
\end{lemma}

As the end of this section, we study the limit behavior of $E_m$ when $m>0$ tends respectively to zero and infinity.
\begin{lemma}\label{lemma:Em_atzero}
  Assume that $N \geq 1$ and $f$ satisfies $(f0)-(f4)$. Then $E_m\to+\infty$ as $m\to 0^+$.
\end{lemma}
\proof It is sufficient to show that for any sequence $\{u_n\}\subset H^1(\mathbb{R}^N)\setminus\{0\}$ such that
  \begin{linenomath*}
    \begin{equation*}
      P(u_n)=0\qquad\text{and}\qquad \lim_{n\to\infty}\|u_n\|_{L^2(\mathbb{R}^N)}=0,
    \end{equation*}
  \end{linenomath*}
one has $I(u_n)\to+\infty$ as $n\to\infty$.  Set
  \begin{linenomath*}
    \begin{equation*}
      s_n:=\ln\bigl(\|\nabla u_n\|_{L^2(\mathbb{R}^N)}\bigr)\qquad\text{and}\qquad v_n:=(-s_n)\star u_n.
    \end{equation*}
  \end{linenomath*}
Clearly, $\|\nabla v_n\|_{L^2(\mathbb{R}^N)}=1$ and $\|v_n\|_{L^2(\mathbb{R}^N)}=\|u_n\|_{L^2(\mathbb{R}^N)}\to0$. Noting that $v_n\to0$ in $L^{2+4/N}(\mathbb{R}^N)$, by Lemma \ref{lemma:I_1} $(ii)$, we have
  \begin{linenomath*}
    \begin{equation*}
      \lim_{n\to\infty}e^{-Ns}\int_{\mathbb{R}^N}F(e^{Ns/2}v_n)dx=0\qquad\text{for any}~s\in\mathbb{R}.
    \end{equation*}
  \end{linenomath*}
Since $P(s_n\star v_n)=P(u_n)=0$, using Lemma \ref{lemma:I_3} $(i)$ and $(ii)$, we derive
  \begin{linenomath*}
    \begin{equation*}
      \begin{split}
        I(u_n)=I(s_n\star v_n)
        &\geq I(s\star v_n)\\
        &=\frac{1}{2}e^{2s}-e^{-Ns}\int_{\mathbb{R}^N}F(e^{Ns/2}v_n)dx=\frac{1}{2}e^{2s}+o_n(1).
      \end{split}
    \end{equation*}
  \end{linenomath*}
As $s\in\mathbb{R}$ is arbitrary, it is clear that $I(u_n)\to+\infty$. \hfill $\square$

\begin{lemma}\label{lemma:Em_atinfinity}
   Assume that $N \geq 1$ and $f$ satisfies $(f6)$ in addition to $(f0)-(f4)$. Then $E_m \to 0$ as $m\to\infty$.
\end{lemma}
\proof  Fix $u \in S_1 \cap L^\infty(\mathbb{R}^N)$ and set $u_m:=\sqrt{m} u \in S_m$ for any $m>1$. By Lemma \ref{lemma:I_3} $(i)$, there exists a unique $s(m) \in \mathbb{R}$ such that $s(m) \star u_m \in \mathcal{P}_m$. Since $F$ is nonnegative by Lemma \ref{lemma:fF}, we then have
  \begin{linenomath*}
    \begin{equation*}
      0<E_m \leq I\big( s(m)\star u_m\big) \leq \frac{1}{2}m e^{2 s(m)}\int_{\mathbb{R}^N}|\nabla u|^2dx.
    \end{equation*}
  \end{linenomath*}
To complete the proof, it is sufficient to show that
  \begin{linenomath*}
    \begin{equation}\label{eq:limit1}
      \lim_{m \to \infty}\sqrt{m} e^{s(m)} = 0.
    \end{equation}
  \end{linenomath*}
Recalling the function $g$ defined by \eqref{eq:g}, from $P(s(m) \star u_m) = 0$, it follows that
  \begin{linenomath*}
    \begin{equation*}
      \int_{\mathbb{R}^N}|\nabla u|^2dx = m^{\frac{2}{N}}\int_{\mathbb{R}^N}g\big( \sqrt{m} e^{N s(m)/2} u\big)|u|^{2+\frac{4}{N}}dx
    \end{equation*}
  \end{linenomath*}
and thus
  \begin{linenomath*}
    \begin{equation}\label{eq:limit2}
      \lim_{m \to \infty}\sqrt{m} e^{N s(m)/2} = 0.
    \end{equation}
  \end{linenomath*}
In particular, we get \eqref{eq:limit1} when $N=1, 2$. Let $\varepsilon >0$ be arbitrarily small. When $N \geq 3$, Lemma \ref{lemma:fF} and $(f6)$ imply that there exists $\delta >0$ small enough such that $\widetilde{F}(t) \geq \frac{4}{N}F(t)\geq\varepsilon ^{-1} |t|^{\frac{2N}{N-2}}$ for any $|t| \leq \delta$. In view of that $P(s(m) \star u_m) = 0$ and \eqref{eq:limit2}, we obtain
  \begin{linenomath*}
    \begin{equation*}
      \begin{split}
        \int_{\mathbb{R}^N}|\nabla u|^2dx
        &= m^{-1}e^{-(N+2)s(m)}\int_{\mathbb{R}^N}\widetilde{F}\big( \sqrt{m} e^{N s(m)/2} u\big)dx\\
        &\geq \varepsilon^{-1} \big[ \sqrt{m} e^{s(m)} \big]^{\frac{4}{N-2}}\int_{\mathbb{R}^N}|u|^{\frac{2N}{N-2}}dx \qquad \text{for large enough}~m
      \end{split}
    \end{equation*}
  \end{linenomath*}
and thus \eqref{eq:limit1} holds when $N \geq 3$. \hfill $\square$
\begin{remark}\label{remark:f6'}
  When $N\geq3$, one still has $\lim_{m\to\infty}E_m=0$ assuming only that $f(t)t/|t|^{\frac{2N}{N-2}} \to +\infty$ as $t\to0^+$ (or as $t\to 0^-$). Indeed, we just need to choose a nonnegative (or nonpositive) function $u \in S_1 \cap L^\infty(\mathbb{R}^N)$ in the proof of Lemma \ref{lemma:Em_atinfinity}.
\end{remark}
\begin{remark}\label{remark:robust}
   When studying $L^2$ constrained mass supercritical problems (set on $\mathbb{R}^N$), the existence of ground states is particularly relevant in view of their physical interpretation.  Assume that we have identified a possible ground state level, say  $E_m$, and managed to find an associated non-vanishing bounded Palais-Smale sequence $\{u_n\}$. Since the working space in general does not embed compactly into any space $L^p(\mathbb{R}^N)$, recovering the compactness of the sequence $\{u_n\}$ may be troublesome. To overcome this difficulty, a by now standard strategy is the one initially proposed in \cite{BJL13}. Roughly speaking, it is to show that $E_m$ is nonincreasing in $m>0$ and satisfies a property similar to Lemma \ref{lemma:Em_strictlymonotonic} and then determine that the Lagrange multiplier is positive. One should however note that, to prove the analogue results to Lemmas \ref{lemma:Em_nonincreasing} and \ref{lemma:Em_strictlymonotonic} (as well as to Lemma \ref{lemma:Em_continuity}), the previous arguments in the literature rely essentially on the homogeneity and the $C^1$ regularity of the nonlinearity $f$ and so the effectiveness of the above strategy had only been confirmed for power type nonlinearities. In the present paper, we develop new arguments for the proofs of Lemmas \ref{lemma:Em_nonincreasing} and \ref{lemma:Em_strictlymonotonic} (as well as of Lemma \ref{lemma:Em_continuity}) which are robust in the sense that they work for the nonlinearity $f$ which is highly non-homogeneous and only continuous.  They should likely be useful to consider other $L^2$ constrained equations in general mass supercritical settings.
\end{remark}

\section{Ground states}\label{sect:groundstate}
In this section we establish the existence of ground states to \eqref{P_m} and complete the study of the properties of the function $m\mapsto E_m$. We deal with the proof of Theorem \ref{theorem:groundstate} at first.  Recall that $N\geq1$, $f$ satisfies $(f0)-(f4)$, and that
  \begin{linenomath*}
    \begin{equation*}
      E_m:=\inf_{u\in\mathcal{P}_m}I(u)>0
    \end{equation*}
  \end{linenomath*}
by Lemma \ref{lemma:I_4} $(iii)$.  As already pointed out in the introduction,  one of the key ingredients for the proof of Theorem \ref{theorem:groundstate} is the following result.
\begin{lemma}\label{lemma:I_5}
  There exists a Palais-Smale sequence $\{u_n\}\subset \mathcal{P}_m$ for the constrained functional $I_{|S_m}$ at the level $E_m$. When $f$ is odd, we have in addition $\|u^-_n\|_{L^2(\mathbb{R}^N)} \to 0$, where $v^-$ stands for the negative part of $v$.
\end{lemma}

To prove Lemma \ref{lemma:I_5}, we borrow some arguments from \cite{BS18,BS19}. Let us first introduce the free functional $\Psi:H^1(\mathbb{R}^N)\setminus\{0\}\to\mathbb{R}$ defined by
  \begin{linenomath*}
    \begin{equation*}\label{eq:Psi}
      \Psi(u):=I(s(u)\star u)=\frac{1}{2}e^{2s(u)}\int_{\mathbb{R}^N}|\nabla u|^2dx-e^{-Ns(u)}\int_{\mathbb{R}^N}F(e^{Ns(u)/2}u)dx,
    \end{equation*}
  \end{linenomath*}
where $s(u)\in\mathbb{R}$ is the unique number guaranteed by Lemma \ref{lemma:I_3}. Inspired by \cite[Proposition 2.9]{SW09} (see also \cite[Proposition 9]{SW10}), we observe
\begin{lemma}\label{lemma:Psi}
  The functional $\Psi:H^1(\mathbb{R}^N)\setminus\{0\}\to\mathbb{R}$ is of class $C^1$ and
    \begin{linenomath*}
      \begin{equation*}
        \begin{split}
          d\Psi(u)[\varphi]&=e^{2s(u)}\int_{\mathbb{R}^N}\nabla u\cdot \nabla \varphi dx-e^{-Ns(u)}\int_{\mathbb{R}^N}f(e^{Ns(u)/2}u)e^{Ns(u)/2}\varphi dx\\
          &=dI(s(u)\star u)[s(u)\star \varphi]
        \end{split}
      \end{equation*}
    \end{linenomath*}
  for any $u\in H^1(\mathbb{R}^N)\setminus\{0\}$ and $\varphi\in H^1(\mathbb{R}^N)$.
\end{lemma}
\proof Let $u\in H^1(\mathbb{R}^N)\setminus\{0\}$ and $\varphi\in H^1(\mathbb{R}^N)$. We estimate the term
  \begin{linenomath*}
    \begin{equation*}
      \Psi(u+t\varphi)-\Psi(u)=I(s_t\star (u+t\varphi))-I(s_0\star u),
    \end{equation*}
  \end{linenomath*}
where $|t|$ is small enough and $s_t:=s(u+t\varphi)$. By the fact that $s_0=s(u)$ is the unique maxmimum point of the function $I(s\star u)$ and the mean value theorem, we have
  \begin{linenomath*}
    \begin{equation*}
      \begin{split}
        I(s_t\star &(u+t\varphi))-I(s_0\star u)
        \leq I(s_t\star (u+t\varphi))-I(s_t\star u)\\
        &=\frac{1}{2}e^{2s_t}\int_{\mathbb{R}^N}\Big[|\nabla (u+t\varphi)|^2-|\nabla u|^2\Big]dx-e^{-Ns_t}\int_{\mathbb{R}^N}\Big[F(e^{Ns_t/2}(u+t\varphi))-F(e^{Ns_t/2}u)\Big]dx\\
        &=\frac{1}{2}e^{2s_t}\int_{\mathbb{R}^N}\bigl(2t\nabla u \cdot \nabla \varphi+t^2|\nabla \varphi|^2\bigr)dx-e^{-Ns_t}\int_{\mathbb{R}^N}f(e^{Ns_t/2}(u+\eta_tt\varphi))e^{Ns_t/2}t\varphi dx,
      \end{split}
    \end{equation*}
  \end{linenomath*}
where $\eta_t\in(0,1)$. Similarly,
  \begin{linenomath*}
    \begin{equation*}
      \begin{split}
        I(s_t\star (u+&t\varphi))-I(s_0\star u)
        \geq I(s_0\star (u+t\varphi))-I(s_0\star u)\\
        &=\frac{1}{2}e^{2s_0}\int_{\mathbb{R}^N}\bigl(2t\nabla u \cdot \nabla \varphi+t^2|\nabla \varphi|^2\bigr)dx-e^{-Ns_0}\int_{\mathbb{R}^N}f(e^{Ns_0/2}(u+\tau_tt\varphi))e^{Ns_0/2}t\varphi dx,
      \end{split}
    \end{equation*}
  \end{linenomath*}
where $\tau_t\in(0,1)$. Since $\lim_{t\to0}s_t= s_0=s(u)$ by Lemma \ref{lemma:I_3} $(iii)$, from the two inequalities above, it follows that
  \begin{linenomath*}
    \begin{equation*}
      \lim_{t\to0}\frac{\Psi(u+t\varphi)-\Psi(u)}{t}
      =e^{2s(u)}\int_{\mathbb{R}^N}\nabla u\cdot \nabla \varphi dx-e^{-Ns(u)}\int_{\mathbb{R}^N}f(e^{Ns(u)/2}u)e^{Ns(u)/2}\varphi dx.
    \end{equation*}
  \end{linenomath*}
By Lemma \ref{lemma:I_3} $(iii)$ again, we see that
the G\^{a}teaux derivative of $\Psi$ is bounded linear in $\varphi$ and continuous in $u$. Therefore $\Psi$ is of class $C^1$, see e.g. \cite{Ch05,Wi96}. In particular, by changing variables in the integrals, we have
  \begin{linenomath*}
    \begin{equation*}
      \begin{split}
        d\Psi(u)[\varphi]
        &=\int_{\mathbb{R}^N}\nabla (s(u)\star u)\cdot \nabla (s(u)\star\varphi) dx-\int_{\mathbb{R}^N}f(s(u)\star u))(s(u)\star\varphi)dx\\
        &=dI(s(u)\star u)[s(u)\star \varphi].
      \end{split}
    \end{equation*}
  \end{linenomath*}
The proof is complete. \hfill $\square$

For given $m>0$, we now consider the constrained functional
  \begin{linenomath*}
    \begin{equation*}
      J:=\Psi_{|S_m}:S_m\to\mathbb{R}
    \end{equation*}
  \end{linenomath*}
which clearly satisfies
\begin{lemma}\label{lemma:J}
  The functional $J: S_m \to \mathbb{R}$ is of class $C^1$ and
    \begin{linenomath*}
      \begin{equation*}
        dJ(u)[\varphi]=d\Psi(u)[\varphi]=dI(s(u)\star u)[s(u)\star \varphi]
      \end{equation*}
    \end{linenomath*}
  for any $u\in S_m$ and $\varphi\in T_uS_m$.
\end{lemma}

We recall below a definition from \cite{G93} and then establish a technical result showing that a ``nice" minimax value of $J$ will yield a Palais-Smale sequence for the constrained functional $I_{|S_m}$ at the same level made of elements of $\mathcal{P}_m$. After that the proof of Lemma \ref{lemma:I_5} will follow.
\begin{definition}[{\cite[Definition 3.1]{G93}}]\label{definition:HSF}
  Let $B$ be a closed subset of a metric space $X$. We say that a class $\mathcal{G}$ of compact subsets of $X$ is a homotopy stable family with closed boundary $B$ provided
    \begin{itemize}
      \item[$(i)$] every set in $\mathcal{G}$ contains $B$,
      \item[$(ii)$] for any set $A\in\mathcal{G}$ and any homotopy $\eta\in C([0,1]\times X, X)$ that satisfies $\eta(t,u)=u$ for all $(t,u)\in \bigl(\{0\}\times X\bigr)\cup \bigr([0,1]\times B\bigr)$, one has $\eta(\{1\}\times A)\in\mathcal{G}$.
    \end{itemize}
  We remark that the case $B=\emptyset$ is admissible.
\end{definition}
\begin{lemma}\label{lemma:J_1}
  Let $\mathcal{G}$ be a homotopy stable family of compact subsets of $S_m$ (with $B=\emptyset$) and set
    \begin{linenomath*}
      \begin{equation*}
        E_{m,\mathcal{G}}:=\inf_{A\in\mathcal{G}}\max_{u\in A}J(u).
      \end{equation*}
    \end{linenomath*}
  If $E_{m,\mathcal{G}}>0$, then there exists a Palais-Smale sequence $\{u_n\}\subset \mathcal{P}_m$ for the constrained functional $I_{|S_m}$ at the level $E_{m,\mathcal{G}}$. In the particular case when $f$ is odd and $\mathcal{G}$ is the class of all singletons included in $S_m$, we have in addition that $\|u^-_n\|_{L^2(\mathbb{R}^N)} \to 0$.
\end{lemma}
\proof  Let $\{A_n\}\subset \mathcal{G}$ be an arbitrary minimizing sequence of $E_{m,\mathcal{G}}$. We define the mapping
  \begin{linenomath*}
    \begin{equation*}
      \eta:[0,1]\times S_m\to S_m,\qquad\eta(t,u)=(t s(u))\star u,
    \end{equation*}
  \end{linenomath*}
which is continuous by Lemma \ref{lemma:I_3} $(iii)$ and satisfies $\eta(t,u)=u$ for all $(t,u)\in \{0\}\times S_m$. Thus, by the definition of $\mathcal{G}$, one has
  \begin{linenomath*}
    \begin{equation}\label{eq:D_n}
      D_n:=\eta(1,A_n)=\{s(u)\star u~|~u\in A_n\}\in \mathcal{G}.
    \end{equation}
  \end{linenomath*}
In particular, $D_n \subset \mathcal{P}_m$ for every $n\in\mathbb{N}^+$. Since $J(s\star u)=J(u)$ for any $s\in\mathbb{R}$ and any $u\in S_m$, it is clear that $\max_{D_n}J=\max_{A_n}J\to E_{m,\mathcal{G}}$ and thus $\{D_n\}\subset \mathcal{G}$ is another minimizing sequence of $E_{m,\mathcal{G}}$. Now, using the minimax principle \cite[Theorem 3.2]{G93}, we obtain a Palais-Smale sequence $\{v_n\}\subset S_m$ for $J$ at the level $E_{m,\mathcal{G}}$ such that $\text{dist}_{H^1(\mathbb{R}^N)}(v_n,D_n)\to0$ as $n\to\infty$. Let
  \begin{linenomath*}
    \begin{equation*}
      s_n:=s(v_n)\qquad \text{and} \qquad u_n:=s_n\star v_n=s(v_n)\star v_n.
    \end{equation*}
  \end{linenomath*}
We prove below a claim concerning $\{e^{-2s_n}\}$ and then show that $\{u_n\}\subset \mathcal{P}_m$ is the desired sequence.

\smallskip

\textbf{Claim.}  \emph{There exists $C>0$ such that $e^{-2s_n}\leq C$ for every $n$.}

We observe that
  \begin{linenomath*}
    \begin{equation*}
      e^{-2s_n}=\frac{\int_{\mathbb{R}^N}|\nabla v_n|^2dx}{\int_{\mathbb{R}^N}|\nabla u_n|^2dx}.
    \end{equation*}
  \end{linenomath*}
Since $\{u_n\}\subset \mathcal{P}_m$, by Lemma \ref{lemma:I_4} $(ii)$, it is clear that $\{\|\nabla u_n\|_{L^2(\mathbb{R}^N)}\}$ is bounded from below by a positive constant. Regarding the term of $\{v_n\}$, since $D_n\subset \mathcal{P}_m$ for every $n$, we have
  \begin{linenomath*}
    \begin{equation*}
      \max_{D_n}I=\max_{D_n}J\to E_{m,\mathcal{G}}
    \end{equation*}
  \end{linenomath*}
and thus $\{D_n\}$ is uniformly bounded in $H^1(\mathbb{R}^N)$ by Lemma \ref{lemma:I_4} $(iv)$; from $\text{dist}_{H^1(\mathbb{R}^N)}(v_n,D_n)\to0$, it then follows that $\sup_{n}\|\nabla v_n\|_{L^2(\mathbb{R}^N)}<\infty$. Clearly, this proves the Claim.

Now, from $\{u_n\}\subset \mathcal{P}_m$, it follows that
  \begin{linenomath*}
    \begin{equation*}
      I(u_n)=J(u_n)=J(v_n)\to E_{m,\mathcal{G}}.
    \end{equation*}
  \end{linenomath*}
We then only need to show that $\{u_n\}$ is a Palais-Smale sequence for $I$ on $S_m$. For any $\psi\in T_{u_n}S_m$, we have
  \begin{linenomath*}
    \begin{equation*}
      \int_{\mathbb{R}^N}v_n[(-s_n)\star\psi]dx=\int_{\mathbb{R}^N}(s_n\star v_n)\psi dx=\int_{\mathbb{R}^N}u_n\psi dx=0,
    \end{equation*}
  \end{linenomath*}
which means $(-s_n)\star \psi\in T_{v_n}S_m$. Also, $\|(-s_n)\star \psi\|_{H^1}\leq \max\{\sqrt{C},1\}\|\psi\|_{H^1}$ by the Claim. Denoting by $\|\cdot\|_{u,*}$ the dual norm of $(T_{u}S_m)^*$ and using Lemma \ref{lemma:J}, we deduce that
  \begin{linenomath*}
    \begin{equation*}
      \begin{split}
        \|dI(u_n)\|_{u_n,*}
        &=\sup_{\psi\in T_{u_n}S_m, \|\psi\|_{H^1}\leq1}\bigl|dI(u_n)[\psi]\bigr|\\
        &=\sup_{\psi\in T_{u_n}S_m, \|\psi\|_{H^1}\leq1}\bigl|dI(s_n\star v_n)[s_n\star ((-s_n)\star \psi)]\bigr|\\
        &=\sup_{\psi\in T_{u_n}S_m, \|\psi\|_{H^1}\leq1}\bigl|dJ(v_n)[(-s_n)\star \psi]\bigr|\\
        &\leq \|dJ(v_n)\|_{v_n,*} \cdot \sup_{\psi\in T_{u_n}S_m, \|\psi\|_{H^1}\leq1}\|(-s_n)\star \psi\|_{H^1}\\
        &\leq \max\{\sqrt{C},1\} \|dJ(v_n)\|_{v_n,*}.
      \end{split}
    \end{equation*}
  \end{linenomath*}
Since $\{v_n\}\subset S_m$ is a Palais-Smale sequence of $J$, it follows that $\|dI(u_n)\|_{u_n,*}\to0$.

Finally, note that the class of all singletons included in $S_m$ is a homotopy stable family of compact subsets of $S_m$ (with $B=\emptyset$). When $f$ is odd, making this particular choice for $\mathcal{G}$ and noting that $J(u)$ is even in $u \in S_m$ by Lemma \ref{lemma:I_3} $(iv)$,  in the above proof we can choose a minimizing sequence $\{A_n\}\subset\mathcal{G}$
which consists of nonnegative functions (rather than an arbitrary one)  and thus the sequence $\{D_n\}$  defined in \eqref{eq:D_n} inherits this property. Since $\text{dist}_{H^1(\mathbb{R}^N)}(v_n,D_n)\to0$, we obtain a Palais-Smale sequence $\{u_n\}\subset \mathcal{P}_m$ for $I_{|S_m}$  at the level $E_{m,\mathcal{G}}$  satisfying the additional property
    \begin{linenomath*}
      \begin{equation*}
        \|u^-_n\|^2_{L^2(\mathbb{R}^N)}=\|s(v_n)\star v^-_n\|^2_{L^2(\mathbb{R}^N)}=\|v^-_n\|^2_{L^2(\mathbb{R}^N)}\to 0.
      \end{equation*}
    \end{linenomath*}
The proof of this lemma is complete.
\hfill $\square$

\medskip
\noindent
\textbf{Proof of Lemma \ref{lemma:I_5}.} We make use of Lemma \ref{lemma:J_1} in the particular case where $\mathcal{G}$ is the class of all singletons included in $S_m$. Since $E_m>0$, it only remains to show that $E_{m,\mathcal{G}}=E_m$. First note that
  \begin{linenomath*}
    \begin{equation*}
      E_{m,\mathcal{G}}=\inf_{A\in \mathcal{G}}\max_{u\in A}J(u)=\inf_{u\in S_m}I(s(u)\star u).
    \end{equation*}
  \end{linenomath*}
For any $u\in S_m$, we deduce from $s(u)\star u\in \mathcal{P}_m$ that $I(s(u)\star u)\geq E_m$. Therefore, $E_{m,\mathcal{G}}\geq E_m$. On the other hand, for any $u\in\mathcal{P}_m$, we have $s(u)=0$ and thus $I(u)=I(0\star u)\geq E_{m,\mathcal{G}}$; clearly, this implies $E_m\geq E_{m,\mathcal{G}}$. \hfill $\square$

\begin{lemma}\label{lemma:decomposition}
  Let $\{u_n\}\subset S_m$ be any bounded Palais-Smale sequence for the constrained functional $I_{|S_m}$ at the level $E_m >0$, satisfying $P(u_n) \to 0$.  Assume in addition that one of following conditions holds:
    \begin{itemize}
      \item[$(i)$] the condition $(f5)$,
      \item[$(ii)$] $\|u^-_n\|_{L^2(\mathbb{R}^N)} \to 0$ and $N=3,4$.
    \end{itemize}
  Then there exists $u \in S_m$ and $\mu>0$ such that, up to the extraction of a subsequence and up to translations in $\mathbb{R}^N$, $u_n \to u$ strongly in $H^1(\mathbb{R}^N)$ and $-\Delta u+\mu u =f(u)$.
\end{lemma}
\proof  Since $\{u_n\}\subset S_m$ is bounded in $H^1(\mathbb{R}^N)$, without loss of generality, one may assume that $\lim_{n\to\infty}\|\nabla u_n\|_{L^2(\mathbb{R}^N)}$, $\lim_{n\to\infty}\int_{\mathbb{R}^N}F(u_n)dx$ and $\lim_{n\to\infty}\int_{\mathbb{R}^N}f(u_n)u_ndx$ exist. Also, from the condition that $\|dI(u_n)\|_{u_n,*}\to0$ and \cite[Lemma 3]{Be83-2}, it follows that
  \begin{linenomath*}
    \begin{equation*}
      -\Delta u_n+\mu_n u_n-f(u_n)\to0\qquad\text{in}~(H^1(\mathbb{R}^N))^*,
    \end{equation*}
  \end{linenomath*}
where
  \begin{linenomath*}
    \begin{equation*}
      \mu_n:=\frac{1}{m}\Big(\int_{\mathbb{R}^N}f(u_n)u_ndx-\int_{\mathbb{R}^N}|\nabla u_n|^2dx\Big).
    \end{equation*}
  \end{linenomath*}
Noting that $\mu_n\to\mu$ for some $\mu\in\mathbb{R}$, we have
  \begin{linenomath*}
    \begin{equation}\label{eq:PS}
      -\Delta u_n(\cdot+y_n)+\mu u_n(\cdot+y_n)-f(u_n(\cdot+y_n))\to0\qquad\text{in}~(H^1(\mathbb{R}^N))^*
    \end{equation}
  \end{linenomath*}
for any $\{y_n\}\subset \mathbb{R}^N$. As a stepping stone, we claim that $\{u_n\}$ is non-vanishing. Indeed, if $\{u_n\}$ is vanishing then $u_n\to0$ in $L^{2+4/N}(\mathbb{R}^N)$ by Lions Lemma \cite[Lemma I.1]{Lions84-2}.  In view of Lemma \ref{lemma:I_1} $(ii)$ and that $P(u_n) \to 0$, we have $\int_{\mathbb{R}^N}F(u_n)dx\to0$ and
  \begin{linenomath*}
    \begin{equation*}
      \int_{\mathbb{R}^N}|\nabla u_n|^2dx=P(u_n)+\frac{N}{2}\int_{\mathbb{R}^N}\widetilde{F}(u_n)dx\to0.
    \end{equation*}
  \end{linenomath*}
As a consequence,
  \begin{linenomath*}
    \begin{equation*}
      E_m=\lim_{n\to\infty}I(u_n)=\frac{1}{2}\lim_{n\to\infty}\int_{\mathbb{R}^N}|\nabla u_n|^2dx-\lim_{n\to\infty}\int_{\mathbb{R}^N}F(u_n)dx=0
    \end{equation*}
  \end{linenomath*}
contradicting the fact that $E_m>0$, and so the claim follows.  The sequence $\{u_n\}$ being non-vanishing, up to a subsequence, there exists $\{y^1_n\}\subset\mathbb{R}^N$ and $w^1\in B_m\setminus\{0\}$
 such that $u_n(\cdot+y^1_n)\rightharpoonup w^1$ in $H^1(\mathbb{R}^N)$, $u_n(\cdot+y^1_n)\to w^1$ in $L^p_{\text{loc}}(\mathbb{R}^N)$ for any $p\in[1,2^*)$, and $u_n(\cdot+y^1_n)\to w^1$ almost everywhere in $\mathbb{R}^N$. Since $f$ satisfies the conditions $(C1)-(C3)$ in Lemma \ref{lemma:BL} by $(f0)-(f2)$, with the aid of \cite[Compactness Lemma 2]{St77} (or \cite[Lemma A.I]{Be83-1}), one has
  \begin{linenomath*}
    \begin{equation*}
      \begin{split}
        \lim_{n\to\infty}\int_{\mathbb{R}^N}\bigl|\bigl[f(u_n&(\cdot+y^1_n))-f(w^1)\bigr]\varphi\bigr|dx\\
        &~~~~\leq\|\varphi\|_{L^\infty(\mathbb{R}^N)}\lim_{n\to\infty}\int_{\text{supp}(\varphi)}\bigl|f(u_n(\cdot+y^1_n))-f(w^1)\bigr|dx=0
      \end{split}
    \end{equation*}
  \end{linenomath*}
for any $\varphi\in C^\infty_0(\mathbb{R}^N)$. In view of \eqref{eq:PS}, we obtain
  \begin{linenomath*}
    \begin{equation}\label{eq:solution}
      -\Delta w^1+\mu w^1=f(w^1).
    \end{equation}
  \end{linenomath*}
In particular, $P(w^1)=0$ by the Nehari and Pohozaev identities corresponding to \eqref{eq:solution}. Let $v^1_n:=u_n-w^1(\cdot-y^1_n)$ for every $n\in\mathbb{N}^+$. Clearly, $v^1_n(\cdot+y^1_n)\rightharpoonup 0$ in $H^1(\mathbb{R}^N)$ and thus
  \begin{linenomath*}
    \begin{equation}\label{eq:L1}
      \begin{split}
        m=\lim_{n\to\infty}\|v^1_n(\cdot+y^1_n)+w^1\|^2_{L^2(\mathbb{R}^N)}
        =\|w^1\|^2_{L^2(\mathbb{R}^N)}+\lim_{n\to\infty}\|v^1_n\|^2_{L^2(\mathbb{R}^N)}.
      \end{split}
    \end{equation}
  \end{linenomath*}
By Lemma \ref{lemma:BL}, we also have
  \begin{linenomath*}
    \begin{equation*}
      \lim_{n\to\infty}\int_{\mathbb{R}^N}F(u_n(\cdot+y^1_n))dx=\int_{\mathbb{R}^N}F(w^1)dx+\lim_{n\to\infty}\int_{\mathbb{R}^N}F(v^1_n(\cdot+y^1_n))dx.
    \end{equation*}
  \end{linenomath*}
It then follows that
  \begin{linenomath*}
    \begin{equation}\label{eq:L2}
      \begin{split}
        E_m=\lim_{n\to\infty}I(u_n)
        &=\lim_{n\to\infty}I(u_n(\cdot+y^1_n))\\
        &=I(w^1)+\lim_{n\to\infty}I(v^1_n(\cdot+y^1_n))=I(w^1)+\lim_{n\to\infty}I(v^1_n).
      \end{split}
    \end{equation}
  \end{linenomath*}
We claim that $\lim_{n\to\infty} I(v^1_n)\geq0$. To see this, we assume by contradiction that $\lim_{n \to \infty}I(v_n^1) <0$. Then $\{v_n^1\}$ is non-vanishing and, up to a subsequence, there exists a sequence  $\{y^{2}_n\}\subset\mathbb{R}^N$ such that
  \begin{linenomath*}
    \begin{equation*}
       \lim_{n\to\infty}\int_{B(y^{2}_n,1)}|v^1_n|^2dx>0.
    \end{equation*}
  \end{linenomath*}
Consequently $|y^2_n-y^1_n|\to\infty$ (since $v^1_n(\cdot+y^1_n)\to0$ in $L^2_\text{loc}(\mathbb{R}^N)$) and, up to a subsequence, $v^1_n(\cdot+y^2_n)\rightharpoonup w^2$ in $H^1(\mathbb{R}^N)$ for some $w^2\in B_m\setminus\{0\}$. Since
  \begin{linenomath*}
    \begin{equation*}
       u_n(\cdot+y^2_n)=v^1_n(\cdot+y^2_n)+ w^1(\cdot - y^1_n + y^2_n) \rightharpoonup w^2\qquad\text{in}~H^1(\mathbb{R}^N),
    \end{equation*}
  \end{linenomath*}
by \eqref{eq:PS} and arguing as above, we deduce that  $P(w^2)=0$ and thus $I(w^2)>0$. Set $v_n^2:= v^1_n - w^2(\cdot -y_n^2)=u_n - \sum ^2_{l=1}w^l(\cdot -y^l_n)$. It is clear that
  \begin{linenomath*}
    \begin{equation*}
      \lim_{n\to\infty}\|\nabla v^2_n\|^2_{L^2(\mathbb{R}^N)}=\lim_{n\to\infty}\|\nabla u_n\|^2_{L^2(\mathbb{R}^N)}-\sum^2_{l=1}\|\nabla w^l\|^2_{L^2(\mathbb{R}^N)}
    \end{equation*}
  \end{linenomath*}
and
  \begin{linenomath*}
    \begin{equation*}
      0> \lim_{n\to\infty}I(v^1_n)=I(w^2)+\lim_{n\to\infty}I(v^2_n)>\lim_{n\to\infty}I(v^2_n).
    \end{equation*}
  \end{linenomath*}
Proceeding this way successively, we obtain an infinite sequence $\{w^k\}\subset B_m\setminus\{0\}$ such that $P(w^k)=0$ and
  \begin{linenomath*}
    \begin{equation*}
      \sum^k_{l=1}\|\nabla w^l\|^2_{L^2(\mathbb{R}^N)} \leq \lim_{n\to\infty}\|\nabla u_n\|^2_{L^2(\mathbb{R}^N)}<+\infty
    \end{equation*}
  \end{linenomath*}
for any $k\in\mathbb{N}^+$. However this is impossible since Remark \ref{remark:P} implies that there exists a $\delta >0$ such that  $\|\nabla w\|_{L^2(\mathbb{R}^N)} \geq \delta $ for any $w \in B_m \backslash \{0\}$ satisfying $P(w)=0$. Therefore, the claim that $\lim_{n\to\infty} I(v^1_n)\geq0$ is proved.

Now we set $s:=\|w^1\|^2_{L^2(\mathbb{R}^N)} \in (0, m]$.  Since $\lim_{n\to\infty} I(v^1_n)\geq0$ and $w^1\in \mathcal{P}_s$, it follows from \eqref{eq:L2} that
  \begin{linenomath*}
    \begin{equation*}
      E_m=I(w^1)+\lim_{n\to\infty}I(v^1_n)\geq I(w^1)\geq E_s.
    \end{equation*}
  \end{linenomath*}
Noting that $E_m$ is nonincreasing in $m>0$ by Lemma \ref{lemma:Em_nonincreasing}, one then has
  \begin{linenomath*}
    \begin{equation}\label{eq:Es=Em}
      I(w^1)=E_s=E_m
    \end{equation}
  \end{linenomath*}
and
  \begin{linenomath*}
    \begin{equation}\label{eq:zero}
      \lim_{n\to\infty} I(v^1_n)=0.
    \end{equation}
  \end{linenomath*}
Clearly, by \eqref{eq:solution}, \eqref{eq:Es=Em} and Lemma \ref{lemma:Em_decreasing}, we derive that $\mu\geq0$. To show that $s=m$, let us prove that $\mu$ is positive.  For clarity, the following two cases are treated separately.

$(i)$ Suppose that $(f5)$ holds. In this case, when $N \geq 3$ the condition $(f5)$  implies that $NF(t)-\frac{N-2}{2}f(t)t>0$ for any $t\neq0$ and when $N=1,2$ this inequality holds thanks to Lemma \ref{lemma:fF} which guarantees that $f(t)t >0$ and $F(t) >0$ for any $t \neq 0$. Thus, from the  Pohozaev identity corresponding to \eqref{eq:solution}, we obtain
  \begin{linenomath*}
    \begin{equation}\label{eq:mu}
      \mu=\frac{1}{m}\int_{\mathbb{R}^N}\Big(NF(w^1)-\frac{N-2}{2}f(w^1)w^1\Big)dx>0.
    \end{equation}
  \end{linenomath*}

$(ii)$ Suppose that $\|u^-_n\|_{L^2(\mathbb{R}^N)} \to 0$ and $N=3,4$. In this case, we assume by contradiction that $\mu =0$. Since $\|u^-_n\|_{L^2(\mathbb{R}^N)} \to 0$ implies $w^1\geq0$,  it follows from \eqref{eq:solution} and Lemma \ref{lemma:fF} that
      \begin{equation*}
        -\Delta w^1=f(w^1)\geq0\quad\text{in}~\mathbb{R}^N.
      \end{equation*}
  Applying \cite[Lemma A.2]{Ik14} with $p:=2 \leq N/(N-2)$, we obtain that $w^1\equiv 0$. This contradicts the fact that $w^1\in S_s$, and thus $\mu>0$.

In both cases, we have proved that $\mu>0$. If $s<m$, taking into account \eqref{eq:solution}, \eqref{eq:Es=Em} and Lemma \ref{lemma:Em_decreasing}, we would have
  \begin{linenomath*}
    \begin{equation*}
      I(w^1)=E_s>E_m
    \end{equation*}
  \end{linenomath*}
which contradicts \eqref{eq:Es=Em}. Therefore, $s:=\|w^1\|^2_{L^2(\mathbb{R}^N)}=m$ and then $\|v^1_n\|_{L^2(\mathbb{R}^N)}\to 0$ via \eqref{eq:L1}. Since $\lim_{n\to\infty}\int_{\mathbb{R}^N}F(v^1_n)dx=0$ by Lemma \ref{lemma:I_1} $(ii)$, we conclude from \eqref{eq:zero} that $\|\nabla v^1_n\|_{L^2(\mathbb{R}^N)}\to 0$ and thus $u_n(\cdot + y_n^1) \to w^1$ strongly in $H^1(\mathbb{R}^N)$. At this point, the proof of the lemma is complete. \hfill $\square$

\begin{remark}\label{remark:mu_positive}
  Showing that $\mu >0$ is crucial to locate the weak limit $w^1$ onto $S_m$ and thus it is one of the key elements to get the strong convergence.  When $N=1,2$, the conditions $(f0)-(f4)$ are sufficient for that purpose; while, when $N \geq 3$, to deal with an arbitrary bounded Palais-Smale sequence satisfying $P(u_n) \to 0$, we need to require in addition $(f5)$.
\end{remark}
\begin{remark}\label{remark:simplification}
  When $f$ satisfies some stronger conditions, it is possible to prove in a simpler way that $\lim_{n\to\infty} I(v^1_n)\geq0$. For example, let us assume that, in addition to $(f0)-(f4)$, $\widetilde{F}$ is of class $C^1$ and that $\widetilde{F}'$ satisfies the conditions $(C1)-(C3)$ in Lemma \ref{lemma:BL}.  Then, as in the proof of \eqref{eq:L2} and noting that $P(w^1)=0$, we have
    \begin{linenomath*}
      \begin{equation*}
        0=\lim_{n\to\infty}P(u_n)=P(w^1)+\lim_{n\to\infty} P(v^1_n)=\lim_{n\to\infty} P(v^1_n)
      \end{equation*}
    \end{linenomath*}
  and thus, using also Lemma \ref{lemma:fF},
    \begin{linenomath*}
      \begin{equation*}
        \lim_{n\to\infty}I(v^1_n)=\lim_{n\to\infty}\Big(I(v^1_n)-\frac{1}{2}P(v^1_n)\Big)=\frac{N}{4}\lim_{n\to\infty}\int_{\mathbb{R}^N}\Big[f(v^1_n)v^1_n-\Big(2+\frac{4}{N}\Big)F(v^1_n)\Big]dx \geq 0.
      \end{equation*}
    \end{linenomath*}
\end{remark}

\medskip

Using Lemmas \ref{lemma:I_5} and \ref{lemma:decomposition}, we are in the position to prove Theorem \ref{theorem:groundstate}.

\noindent
\textbf{Proof of Theorem \ref{theorem:groundstate}.}
By Lemmas \ref{lemma:I_5} and \ref{lemma:I_4} $(iv)$, we have a bounded Palais-Smale sequence $\{u_n\}\subset \mathcal{P}_m$ for the constrained functional $I_{|S_m}$ at the level $E_m >0$.

$(i)$ Suppose that $(f5)$ holds. Then Lemma \ref{lemma:decomposition} applies and provides the existence of a ground state $u \in S_m$ at the level $E_m$.

$(ii)$ Assume that $f$ is odd and that $(f5)$ holds for $N\geq 5$. Then, by Lemma \ref{lemma:I_5}, we have in addition $\|u^-_n\|_{L^2(\mathbb{R}^N)} \to 0$. Applying Lemma \ref{lemma:decomposition}, we obtain a nonnegative ground state $u \in S_m$ at the level $E_m$. Moreover, by the strong maximum principle, $u>0$ as required. \hfill $\square$

\medskip

With Theorem \ref{theorem:groundstate} at hand, we now have all the elements to prove Theorem \ref{theorem:Em}.

\noindent
\textbf{Proof of Theorem \ref{theorem:Em}.} We first prove the strict decrease of $E_m$. In each of the three cases, from Theorem \ref{theorem:groundstate} we know that $E_m$ is reached by a ground state of \eqref{P_m} with the associated Lagrange multiplier being positive, and thus the function $m \mapsto E_m$ is strictly decreasing on $(0, \infty)$ by Lemma \ref{lemma:Em_decreasing}. The rest of the proof directly follows from Lemmas \ref{lemma:I_4}, \ref{lemma:Em_continuity}, \ref{lemma:Em_nonincreasing}, \ref{lemma:Em_atzero} and \ref{lemma:Em_atinfinity}. \hfill $\square$

\begin{remark}\label{remark:f7}
Let us show that, when $N\geq2$ and $f$ satisfies the assumptions of Theorem \ref{theorem:groundstate} $(ii)$ and
\begin{itemize}
  \item[$(f7)$] $h(t):=[f(t)t-(2+4/N)F(t)]/t^2$ is nonincreasing on $(-\infty,0)$ and nondecreasing on $(0,\infty)$,
\end{itemize}
we can work directly in the subspace of radially symmetric functions to derive a positive radial ground state of \eqref{P_m}. Indeed, for given $u\in S_m$ and any $s\in\mathbb{R}$, one has
  \begin{linenomath*}
    \begin{equation*}
      \begin{split}
        I(s\star u)-\frac{1}{2}P(s\star u)
        &=\frac{N}{4}\int_{\mathbb{R}^N}\Big[f(s\star u)s\star u-(2+4/N)F(s\star u)\Big]dx\\
        &=\frac{N}{4}\int_{\mathbb{R}^N}e^{-Ns}\Big[f(e^{Ns/2}u)e^{Ns/2}u-(2+4/N)F(e^{Ns/2}u)\Big]dx\\
        &=\frac{N}{4}\int_{\mathbb{R}^N}\frac{f(e^{Ns/2}u)e^{Ns/2}u-(2+4/N)F(e^{Ns/2}u)}{(e^{Ns/2}u)^2}u^2dx
      \end{split}
    \end{equation*}
  \end{linenomath*}
and then, by $(f7)$, the function $I(s\star u)-\frac{1}{2}P(s\star u)$ is nondecreasing in $s\in\mathbb{R}$. Let
\begin{equation*}
  E_{m,r}:=\inf_{u \in \mathcal{P}_m \cap H^1_r(\mathbb{R}^N)} I(u).
\end{equation*}
For any given $u \in \mathcal{P}_m$, define $\tilde{u}:= |u|^*$ as the Schwarz symmetrization of $|u|$. Clearly, $\tilde{u} \in S_m \cap H^1_r(\mathbb{R}^N)$ and $P(\tilde{u}) \leq P(u) = 0$. By Lemma \ref{lemma:I_3}, there exists $s:= s(\tilde{u}) \leq 0$ such that $P(s\star\tilde{u})=0$ and thus
\begin{linenomath*}
    \begin{equation*}
      \begin{split}
        E_m \leq E_{m,r}\leq I(s\star \tilde{u})
        & = I(s \star \tilde{u}) - \frac{1}{2}P(s \star \tilde{u}) \\
        &\leq I(\tilde{u}) - \frac{1}{2}P(\tilde{u}) \\
        &=\frac{N}{4}\int_{\mathbb{R}^N}\Big(\widetilde{F}(\tilde{u})-\frac{4}{N}F(\tilde{u})\Big)dx\\
        &=\frac{N}{4}\int_{\mathbb{R}^N}\Big(\widetilde{F}(u)-\frac{4}{N}F(u)\Big)dx=I(u).
      \end{split}
    \end{equation*}
  \end{linenomath*}
Since $u \in \mathcal{P}_m$ is arbitrary, it is clear that
  \begin{equation*}
    E_{m,r}=E_m>0.
  \end{equation*}
Now, similarly to the proof of Lemma \ref{lemma:I_5}, we can find a Palais-Smale sequence $\{u_n\}\subset \mathcal{P}_m \cap H^1_r(\mathbb{R}^N)$ for the constrained functional $I_{|S_m \cap H^1_r(\mathbb{R}^N)}$, at the level $E_{m,r}=E_m >0$, satisfying $\|u^-_n\|_{L^2(\mathbb{R}^N)} \to 0$. In view of Lemma \ref{lemma:I_4} $(iv)$, the sequence $\{u_n\}$ is bounded in $H^1_r(\mathbb{R}^N)$. When $N=2$ or $N\geq5$, by a later compactness result Lemma \ref{lemma:compact}, we obtain a positive radial ground state $u \in S_m$ at the level $E_m$. It is notable that the proof of this case does not use Lemmas \ref{lemma:Em_nonincreasing} and \ref{lemma:Em_decreasing}. When $N=3,4$, we can conclude the proof by adapting the argument of Lemma \ref{lemma:decomposition} $(ii)$. In particular, since the inclusion $H^1_r(\mathbb{R}^N)\hookrightarrow L^{2+4/N}(\mathbb{R}^N)$ is compact, $\{y^1_n\}$ is now chosen as the zero sequence and the claim that $\lim_{n\to\infty} I(v^1_n)\geq0$ can be proved more easily.

One should also note that, when $\widetilde{F}(t):=f(t)t-2F(t)$ is of class $C^1$, the condition $(f4)$ implies $(f7)$. Indeed, recalling that $g(t):=\widetilde{F}(t)/|t|^{2+4/N}$, we have
  $$g'(t)=\frac{\widetilde{F}'(t)t-(2+4/N)\widetilde{F}(t)}{|t|^{3+4/N}\text{sign}(t)},$$
and thus, by $(f4)$,
  \begin{equation*}
    \widetilde{F}'(t)t-(2+4/N)\widetilde{F}(t)\geq0\qquad\text{for any}~t\neq0.
  \end{equation*}
Since
  $$h'(t)=\left[\frac{\widetilde{F}(t)-\frac{4}{N}F(t)}{t^2}\right]'=\frac{\widetilde{F}'(t)t-(2+4/N)\widetilde{F}(t)}{t^3},$$
the condition $(f7)$ follows.
\end{remark}

\section{Radial solutions}\label{sect:radial}
This section concerns the existence of infinitely many radial solutions to \eqref{P_m}  when $N\geq2$ and $f$ is an odd function satisfying $(f0)-(f5)$.  To prove Theorem \ref{theorem:radial}, we prepare below several technical lemmas. Denote by $\sigma: H^1(\mathbb{R}^N)\to H^1(\mathbb{R}^N)$ the transformation $\sigma (u)=-u$ and let $X\subset H^1(\mathbb{R}^N)$. A set $A\subset X$ is called $\sigma$-invariant if $\sigma(A)=A$. A homotopy $\eta:[0,1]\times X\to X$ is $\sigma$-equivariant if $\eta (t,\sigma(u))=\sigma(\eta(t,u))$ for all $(t,u)\in [0,1]\times X$. The following definition is \cite[Definition 7.1]{G93}.
\begin{definition}\label{definition:sigmaHSF}
  Let $B$ be a closed $\sigma$-invariant subset of $X\subset H^1(\mathbb{R}^N)$. A class $\mathcal{G}$ of compact subsets of $X$ is said to be a $\sigma$-homotopy stable family with closed boundary $B$ provided
    \begin{itemize}
      \item[$(i)$] every set in $\mathcal{G}$ is $\sigma$-invariant,
      \item[$(ii)$] every set in $\mathcal{G}$ contains $B$,
      \item[$(iii)$] for any set $A\in\mathcal{G}$ and any $\sigma$-equivariant homotopy $\eta\in C([0,1]\times X, X)$ that satisfies $\eta(t,u)=u$ for all $(t,u)\in \bigl(\{0\}\times X\bigr)\cup \bigr([0,1]\times B\bigr)$, one has $\eta(\{1\}\times A)\in\mathcal{G}$.
    \end{itemize}
\end{definition}

Since $f$ is odd, we see from Lemma \ref{lemma:I_3} $(iv)$ that $s(-u)=s(u)$, and thus  the constrained functional
  \begin{linenomath*}
    \begin{equation*}
      J(u)=I(s(u)\star u)=\frac{1}{2}e^{2s(u)}\int_{\mathbb{R}^N}|\nabla u|^2dx-e^{-Ns(u)}\int_{\mathbb{R}^N}F(e^{Ns(u)/2}u)dx
    \end{equation*}
  \end{linenomath*}
is even in $u\in S_m$. That is, $J$ is a $\sigma$-invariant functional on $S_m$. Recall that $H^1_r(\mathbb{R}^N)$ stands for the space of radially symmetric functions in $H^1(\mathbb{R}^N)$. To get the particular Palais-Smale sequences of $I_{|S_m\cap H^1_r(\mathbb{R}^N)}$ that consist of elements in $\mathcal{P}_m$, we need
\begin{lemma}\label{lemma:J_2}
  Let $\mathcal{G}$ be a $\sigma$-homotopy stable family of compact subsets of $S_m\cap H^1_r(\mathbb{R}^N)$ (with $B=\emptyset$) and set
    \begin{linenomath*}
      \begin{equation*}
        E_{m,\mathcal{G}}:=\inf_{A\in\mathcal{G}}\max_{u\in A}J(u).
      \end{equation*}
    \end{linenomath*}
  If $E_{m,\mathcal{G}}>0$, then there exists a Palais-Smale sequence $\{u_n\}\subset \mathcal{P}_m\cap H^1_r(\mathbb{R}^N)$ for the constrained functional $I_{|S_m\cap H^1_r(\mathbb{R}^N)}$ at the level $E_{m,\mathcal{G}}$.
\end{lemma}
\proof This result is an equivariant version of Lemma \ref{lemma:J_1}. The proof is almost identical to that of Lemma \ref{lemma:J_1} but makes use of \cite[Theorem 7.2]{G93} instead of \cite[Theorem 3.2]{G93}. \hfill $\square$

We now construct a sequence of $\sigma$-homotopy stable families of compact subsets of $S_m\cap H^1_r(\mathbb{R}^N)$ (with $B=\emptyset$). Fix a sequence of finite dimensional linear subspaces $\{V_k\}\subset H^1_r(\mathbb{R}^N)$ such that $V_k\subset V_{k+1}$, $\dim V_k=k$ and $\cup_{k\geq1}V_k$ is dense in $H^1_r(\mathbb{R}^N)$, and denote by $\pi_k$ the orthogonal projection from $H^1_r(\mathbb{R}^N)$ onto $V_k$. We also recall the definition of the genus of $\sigma$-invariant sets due to M. A. Krasnoselskii and refer to \cite[Section 7]{Ra86} for its basic properties.
\begin{definition}
  For any nonempty closed $\sigma$-invariant set $A\subset H^1(\mathbb{R}^N)$, the genus of $A$ is defined by
    \begin{linenomath*}
      \begin{equation*}
        \text{Ind}(A):=\min\bigl\{k\in\mathbb{N}^+~|~\exists~\phi:A\to\mathbb{R}^k\setminus\{0\},~\phi~\text{is odd and continuous}\bigr\}.
      \end{equation*}
    \end{linenomath*}
  We set $\text{Ind}(A)=\infty$ if such $\phi$ does not exist, and set $\text{Ind}(A)=0$ if $A=\emptyset$.
\end{definition}
Let $\Sigma$ be the family of compact $\sigma$-invariant subsets of $S_m\cap H^1_r(\mathbb{R}^N)$. For each $k\in\mathbb{N}^+$, set
  \begin{linenomath*}
    \begin{equation*}
      \mathcal{G}_k:=\bigl\{A\in\Sigma~|~\text{Ind}(A)\geq k\bigr\}
    \end{equation*}
  \end{linenomath*}
and
  \begin{linenomath*}
    \begin{equation*}
      E_{m,k}:=\inf_{A\in\mathcal{G}_k}\max_{u\in A}J(u).
    \end{equation*}
  \end{linenomath*}
Concerning $\mathcal{G}_k$ and $E_{m,k}$, we have
\begin{lemma}\label{lemma:GkEmk}
  \begin{itemize}
    \item[$(i)$] For any $k\in\mathbb{N}^+$,
                   \begin{linenomath*}
                     \begin{equation*}
                       \mathcal{G}_k\neq\emptyset
                     \end{equation*}
                   \end{linenomath*}
                 and $\mathcal{G}_k$ is a $\sigma$-homotopy stable family of compact subsets of $S_m\cap H^1_r(\mathbb{R}^N)$ (with $B=\emptyset$).
    \item[$(ii)$] $E_{m,k+1}\geq E_{m,k}>0$ for any $k\in\mathbb{N}^+$.
  \end{itemize}
\end{lemma}
\proof $(i)$ For each $k\in\mathbb{N}^+$, $S_m\cap V_k\in \Sigma$. By the basic properties of the genus, one has
  \begin{linenomath*}
    \begin{equation*}
      \text{Ind}(S_m\cap V_k)=k
    \end{equation*}
  \end{linenomath*}
and thus $\mathcal{G}_k\neq\emptyset$. The rest is clear by Definition \ref{definition:sigmaHSF} and by again the basic properties of the genus.

$(ii)$ By Item $(i)$, $E_{m,k}$ is well-defined. For any $A\in\mathcal{G}_k$, using the fact that $s(u)\star u\in \mathcal{P}_m$ for all $u\in A$ and Lemma \ref{lemma:I_4} $(iii)$, we have
  \begin{linenomath*}
    \begin{equation*}
      \max_{u\in A}J(u)=\max_{u\in A}I(s(u)\star u)\geq \inf_{v\in\mathcal{P}_m}I(v)>0
    \end{equation*}
  \end{linenomath*}
and thus $E_{m,k}>0$. Since $\mathcal{G}_{k+1}\subset\mathcal{G}_k$, it is clear that $E_{m,k+1}\geq E_{m,k}$. \hfill $\square$

Helped by the property that the embedding $H^1_r(\mathbb{R}^N)\hookrightarrow L^p(\mathbb{R}^N)$ is compact for any $2<p<2^*$, we establish below a compactness result.
\begin{lemma}\label{lemma:compact}
  Let $\{u_n\}\subset S_m \cap H^1_r(\mathbb{R}^N)$ be any bounded Palais-Smale sequence for the constrained functional $I_{|S_m \cap H^1_r(\mathbb{R}^N)}$, at an arbitrary level $c>0$, satisfying $P(u_n) \to 0$. Then there exists $u \in S_m\cap H^1_r(\mathbb{R}^N)$ and $\mu>0$ such that, up to the extraction of a subsequence, $u_n \to u$ strongly in $H^1(\mathbb{R}^N)$ and $-\Delta u+\mu u =f(u)$.
\end{lemma}
\proof  Since the sequence is bounded in $H^1_r(\mathbb{R}^N)$, up to a subsequence, there exists $u\in H^1_r(\mathbb{R}^N)$ such that $u_n\rightharpoonup u$ in $H^1_r(\mathbb{R}^N)$, $u_n\to u$ in $L^p(\mathbb{R}^N)$ for any $p\in (2,2^*)$, and $u_n\to u$ almost everywhere in $\mathbb{R}^N$. Also, from $\|dI(u_n)\|_{u_n,*}\to0$ and \cite[Lemma 3]{Be83-2}, it follows that
  \begin{linenomath*}
    \begin{equation}\label{eq:PS_radial}
      -\Delta u_n+\mu_n u_n-f(u_n)\to0\qquad\text{in}~(H^1_r(\mathbb{R}^N))^*,
    \end{equation}
  \end{linenomath*}
where
  \begin{linenomath*}
    \begin{equation*}
      \mu_n:=\frac{1}{m}\Big(\int_{\mathbb{R}^N}f(u_n)u_ndx-\int_{\mathbb{R}^N}|\nabla u_n|^2dx\Big).
    \end{equation*}
  \end{linenomath*}
Without loss of generality, one may assume that $\mu_n\to\mu$ for some $\mu\in\mathbb{R}$. Similarly to the proof of \eqref{eq:solution} and using the Palais principle of symmetric criticality \cite{Pa79}, we obtain
  \begin{linenomath*}
    \begin{equation}\label{eq:solution_radial}
      -\Delta u+\mu u=f(u).
    \end{equation}
  \end{linenomath*}
To proceed further, we claim that $u\neq0$. Indeed, if $u=0$ then $u_n\to0$ in $L^{2+4/N}(\mathbb{R}^N)$. In view of Lemma \ref{lemma:I_1} $(ii)$ and that $P(u_n) \to 0$, we have $\int_{\mathbb{R}^N}F(u_n)dx\to0$ and
  \begin{linenomath*}
    \begin{equation*}
      \int_{\mathbb{R}^N}|\nabla u_n|^2dx=P(u_n)+\frac{N}{2}\int_{\mathbb{R}^N}\widetilde{F}(u_n)dx\to0.
    \end{equation*}
  \end{linenomath*}
As a consequence,
  \begin{linenomath*}
    \begin{equation*}
      c=\lim_{n\to\infty}I(u_n)=\frac{1}{2}\lim_{n\to\infty}\int_{\mathbb{R}^N}|\nabla u_n|^2dx-\lim_{n\to\infty}\int_{\mathbb{R}^N}F(u_n)dx=0
    \end{equation*}
  \end{linenomath*}
which contradicts the condition that $c>0$.  Now, by the fact that $u\neq0$ and similarly to the proof of \eqref{eq:mu}, it is clear that
  \begin{linenomath*}
    \begin{equation*}
      \mu=\frac{1}{m}\int_{\mathbb{R}^N}\Big(NF(u)-\frac{N-2}{2}f(u)u\Big)dx>0.
    \end{equation*}
  \end{linenomath*}
Since $u_n\rightharpoonup u$ in $H^1_r(\mathbb{R}^N)$, one can show in a standard way that $\int_{\mathbb{R}^N}\big[f(u_n)-f(u)\big]udx\to 0$. Noting that $u_n\to u$ in $L^{2+4/N}(\mathbb{R}^N)$, we have $\int_{\mathbb{R}^N}f(u_n)(u_n-u)dx\to0$ by Lemma \ref{lemma:I_1} $(iii)$ and thus
  \begin{linenomath*}
    \begin{equation*}
      \lim_{n\to\infty}\int_{\mathbb{R}^N}f(u_n)u_ndx=\int_{\mathbb{R}^N}f(u)udx.
    \end{equation*}
  \end{linenomath*}
In view of \eqref{eq:solution_radial} and \eqref{eq:PS_radial}, it follows that
  \begin{linenomath*}
    \begin{equation*}
      \begin{split}
         \int_{\mathbb{R}^N}|\nabla u|^2dx+\mu\int_{\mathbb{R}^N}u^2dx
         &=\int_{\mathbb{R}^N}f(u)udx\\
         &=\lim_{n\to\infty}\int_{\mathbb{R}^N}f(u_n)u_ndx=\lim_{n\to\infty}\int_{\mathbb{R}^N}|\nabla u_n|^2dx+\mu m.
      \end{split}
    \end{equation*}
  \end{linenomath*}
Since $\mu>0$, we obtain
  \begin{linenomath*}
    \begin{equation*}
      \lim_{n\to\infty}\int_{\mathbb{R}^N}|\nabla u_n|^2dx=\int_{\mathbb{R}^N}|\nabla u|^2dx,\qquad\lim_{n\to\infty}\int_{\mathbb{R}^N}u_n^2dx=m=\int_{\mathbb{R}^N}u^2dx,
    \end{equation*}
  \end{linenomath*}
and thus $u_n\to u$ in $H^1_r(\mathbb{R}^N)$. \hfill $\square$

The next result concerns the limit behaviour of $E_{m,k}$ when $k\to\infty$ and it serves as an essential and final preparation for the proof of Theorem \ref{theorem:radial}.
\begin{lemma}\label{lemma:Emk}
  $E_{m,k}\to+\infty$ as $k\to\infty$.
\end{lemma}

Since we do not require that $\widetilde{F}$ is of class $C^1$, the Pohozaev manifold $\mathcal{P}_m$ is in general only a topological manifold. Despite the fact that we have Lemma \ref{lemma:compact} and that the constrained functional $I_{|\mathcal{P}_m}$ is bounded from below and coercive by Lemma \ref{lemma:I_4}, it is problematic to prove Lemma \ref{lemma:Emk} by a standard genus type argument for $I_{|\mathcal{P}_m}$.  Our proof of Lemma \ref{lemma:Emk} is inspired by that of \cite[Theorem 9]{Be83-2} and relies on the following Lemma \ref{lemma:unbounded}.
\begin{lemma}\label{lemma:unbounded}
  For any $c>0$, there exists $\rho=\rho(c)>0$ small enough and $k(c)\in \mathbb{N}^+$ sufficiently large such that for any $k\geq k(c)$ and any $u\in\mathcal{P}_m\cap H^1_r(\mathbb{R}^N)$ one has
    \begin{linenomath*}
      \begin{equation*}
        I(u)\geq c\qquad\text{if}~\|\pi_ku\|_{H^1(\mathbb{R}^N)}\leq\rho.
      \end{equation*}
    \end{linenomath*}
\end{lemma}
\proof By contradiction, we assume that there exists $c_0>0$ such that for any $\rho>0$ and any $k\in\mathbb{N}^+$ one can always find $l=l(\rho,k)\geq k$ and $u=u(\rho,k)\in \mathcal{P}_m\cap H^1_r(\mathbb{R}^N)$ such that
  \begin{linenomath*}
    \begin{equation*}
      \|\pi_lu\|_{H^1(\mathbb{R}^N)}\leq\rho\qquad\text{but}\qquad I(u)< c_0.
    \end{equation*}
  \end{linenomath*}
As a consequence, one can obtain a strictly increasing sequence $\{k_j\}\subset\mathbb{N}^+$ (and so $\lim_{j\to\infty}k_j=\infty$) and a sequence $\{u_j\}\subset \mathcal{P}_m\cap H^1_r(\mathbb{R}^N)$ such that
  \begin{linenomath*}
    \begin{equation*}
      \|\pi_{k_j}u_j\|_{H^1(\mathbb{R}^N)}\leq\frac{1}{j}\qquad\text{and}\qquad I(u_j)< c_0
    \end{equation*}
  \end{linenomath*}
for any $j\in\mathbb{N}^+$. Since $\{u_j\}$ is bounded in $H^1_r(\mathbb{R}^N)$ by Lemma \ref{lemma:I_4} $(iv)$, up to a subsequence, there exists $u\in H^1_r(\mathbb{R}^N)$ such that
  \begin{linenomath*}
    \begin{equation*}
      u_j\rightharpoonup u~\text{in}~H^1_r(\mathbb{R}^N)\qquad\text{and}\qquad u_j\rightharpoonup u~\text{in}~L^2(\mathbb{R}^N).
    \end{equation*}
  \end{linenomath*}
To derive a contradiction, we claim that $u=0$. Indeed, from $k_j\to\infty$, it follows that $\pi_{k_j}u\to u$ in $L^2(\mathbb{R}^N)$ and thus
  \begin{linenomath*}
    \begin{equation*}
      \bigl(\pi_{k_j}u_j, u\bigr)_{L^2(\mathbb{R}^N)}=\bigl(u_j, \pi_{k_j}u\bigr)_{L^2(\mathbb{R}^N)}\to\bigl(u, u\bigr)_{L^2(\mathbb{R}^N)}\qquad\text{as}~j\to\infty.
    \end{equation*}
  \end{linenomath*}
Combining the fact that $\pi_{k_j}u_j\to 0$ in $L^2(\mathbb{R}^N)$, we then have
  \begin{linenomath*}
    \begin{equation*}
      \|u\|^2_{L^2(\mathbb{R}^N)}=\lim_{j\to\infty}\bigl(\pi_{k_j}u_j,u\bigr)_{L^2(\mathbb{R}^N)}=0,
    \end{equation*}
  \end{linenomath*}
which proves the claim.  Now, up to a subsequence,  $\|u_j\|_{L^{2+4/N}(\mathbb{R}^N)}\to0$ by the compact inclusion $H^1_r(\mathbb{R}^N)\hookrightarrow L^{2+4/N}(\mathbb{R}^N)$.  Using that $\{u_j\}\subset \mathcal{P}_m\cap H^1_r(\mathbb{R}^N)$ and Lemma \ref{lemma:I_1} $(ii)$, we obtain
  \begin{linenomath*}
    \begin{equation*}
      \int_{\mathbb{R}^N}|\nabla u_j|^2dx=\frac{N}{2}\int_{\mathbb{R}^N}\widetilde{F}(u_j)dx\to0\qquad\text{as}~j\to\infty,
    \end{equation*}
  \end{linenomath*}
which contradicts Lemma \ref{lemma:I_4} $(ii)$. The proof of Lemma \ref{lemma:unbounded} is complete. \hfill $\square$

\medskip
\noindent
\textbf{Proof of Lemma \ref{lemma:Emk}.}  By contradiction, we assume that
  \begin{linenomath*}
    \begin{equation}\label{eq:bounded}
      \liminf_{k\to\infty}E_{m,k}< c\qquad\text{for some}~c>0.
    \end{equation}
  \end{linenomath*}
Let $\rho(c)>0$ and $k(c)\in\mathbb{N}^+$ be the numbers given by Lemma \ref{lemma:unbounded}. Clearly, in view of \eqref{eq:bounded}, there exists $k>k(c)$ such that $E_{m,k}<c$. By the definition of $E_{m,k}$, one can then find $A\in \mathcal{G}_k$ (that is $A\in\Sigma$ and $\text{Ind}(A)\geq k$) such that
  \begin{linenomath*}
    \begin{equation*}
      \max_{u\in A}I(s(u)\star u)=\max_{u\in A}J(u)< c.
    \end{equation*}
  \end{linenomath*}
Since Lemma \ref{lemma:I_3} $(iii)$ and $(iv)$ imply that the mapping $\varphi:A\to \mathcal{P}_m\cap H^1_r(\mathbb{R}^N)$ defined by $\varphi(u)=s(u)\star u$ is odd and continuous, we have $\overline{A}:=\varphi(A)\subset \mathcal{P}_m\cap H^1_r(\mathbb{R}^N)$, $\max_{v\in \overline{A}}I(v)< c$ and
  \begin{linenomath*}
    \begin{equation}\label{eq:overlineA}
      \text{Ind}(\overline{A})\geq \text{Ind}(A)\geq k>k(c).
    \end{equation}
  \end{linenomath*}
Also, from Lemma \ref{lemma:unbounded}, it follows that $\inf_{v\in\overline{A}}\|\pi_{k(c)}v\|_{H^1}\geq \rho(c)>0$. Setting
  \begin{linenomath*}
    \begin{equation*}
      \psi(v)=\frac{1}{\|\pi_{k(c)}v\|_{H^1}}\pi_{k(c)}v\qquad\text{for any}~v\in\overline{A},
    \end{equation*}
  \end{linenomath*}
we obtain an odd continuous mapping $\psi:\overline{A}\to\psi(\overline{A})\subset V_{k(c)}\setminus\{0\}$ and thus
  \begin{linenomath*}
    \begin{equation*}
      \text{Ind}(\overline{A})\leq \text{Ind}(\psi(\overline{A}))\leq k(c)
    \end{equation*}
  \end{linenomath*}
which contradicts \eqref{eq:overlineA}. Therefore, we have $E_{m,k}\to+\infty$ as $k\to\infty$. \hfill $\square$

With all the technical lemmas in place, we can now prove Theorem \ref{theorem:radial}.

\medskip
\noindent
\textbf{Proof of Theorem \ref{theorem:radial}.} For each $k\in\mathbb{N}^+$, by Lemmas \ref{lemma:J_2} and \ref{lemma:GkEmk}, one can find a Palais-Smale sequence $\{u^k_n\}^\infty_{n=1}\subset \mathcal{P}_m\cap H^1_r(\mathbb{R}^N)$ of the constrained functional $I_{|S_m\cap H^1_r(\mathbb{R}^N)}$ at the level $E_{m,k}>0$. By Lemma \ref{lemma:I_4} $(iv)$, the sequence is bounded in $H^1_r(\mathbb{R}^N)$ and thus in view of Lemma \ref{lemma:compact}, we deduce that \eqref{P_m} has a radial solution $u_k$ with $I(u_k)=E_{m,k}$. Also, from Lemma \ref{lemma:GkEmk} $(ii)$ and Lemma \ref{lemma:Emk}, it follows that
  \begin{linenomath*}
    \begin{equation*}
      I(u_{k+1})\geq I(u_k)>0\qquad\text{for any}~k\geq1
    \end{equation*}
  \end{linenomath*}
and $I(u_k)\to+\infty$. \hfill $\square$

\section{Nonradial sign-changing solutions}\label{sect:nonradial}
In this section we focus on nonradial sign-changing solutions of \eqref{P_m} when $N\geq4$ and prove Theorems \ref{theorem:nonradial1} and \ref{theorem:nonradial2}. Since the arguments are similar to those for Theorems \ref{theorem:groundstate} and \ref{theorem:radial}, we just outline the proofs.

\subsection{Proof of Theorem \ref{theorem:nonradial2}}\label{subsect:nonradial2}
Recall that $N=4$ or $N\geq6$, $N-2M\neq1$, $X_2:=H^1_{\mathcal{O}_2}\cap X_\omega$ and $f$ is odd satisfying $(f0)-(f5)$.  Similarly to the proofs of Lemmas \ref{lemma:J_1} and \ref{lemma:J_2}, we have
\begin{lemma}\label{lemma:J_3}
  Let $\mathcal{G}$ be a $\sigma$-homotopy stable family of compact subsets of $S_m\cap X_2$ (with $B=\emptyset$) and set
    \begin{linenomath*}
      \begin{equation*}
        E_{m,\mathcal{G}}:=\inf_{A\in\mathcal{G}}\max_{v\in A}J(v).
      \end{equation*}
    \end{linenomath*}
  If $E_{m,\mathcal{G}}>0$, then there exists a Palais-Smale sequence $\{v_n\}\subset \mathcal{P}_m\cap X_2$ for the constrained functional $I_{|S_m\cap X_2}$ at the level $E_{m,\mathcal{G}}$.
\end{lemma}

Let $\overline{\Sigma}$ be the family of compact $\sigma$-invariant subsets of $S_m\cap X_2$. For each $k\in\mathbb{N}^+$, we set
  \begin{linenomath*}
    \begin{equation*}
      \overline{\mathcal{G}}_k:=\bigl\{A\in\overline{\Sigma}~|~\text{Ind}(A)\geq k\bigr\}
    \end{equation*}
  \end{linenomath*}
and
  \begin{linenomath*}
    \begin{equation}\label{eq:Emk_1}
      \overline{E}_{m,k}:=\inf_{A\in\overline{\mathcal{G}}_k}\max_{v\in A}J(v).
    \end{equation}
  \end{linenomath*}
It is clear that $\overline{\mathcal{G}}_k$ and $\overline{E}_{m,k}$ satisfy
\begin{lemma}\label{lemma:GkEmk_1}
  \begin{itemize}
    \item[$(i)$] For any $k\in\mathbb{N}^+$,
                   \begin{linenomath*}
                     \begin{equation*}
                       \overline{\mathcal{G}}_k\neq\emptyset
                     \end{equation*}
                   \end{linenomath*}
                 and $\overline{\mathcal{G}}_k$ is a $\sigma$-homotopy stable family of compact subsets of $S_m\cap X_2$ (with $B=\emptyset$).
    \item[$(ii)$] $\overline{E}_{m,k+1}\geq \overline{E}_{m,k}>0$ for any $k\in\mathbb{N}^+$.
  \end{itemize}
\end{lemma}

Since the inclusion $X_2\hookrightarrow L^p(\mathbb{R}^N)$ is compact for any $2<p<2N/(N-2)$, see \cite{Lions82} or \cite[Corollary 1.25]{Wi96}, we have the following compactness result by adapting the proof of Lemma \ref{lemma:compact}.
\begin{lemma}\label{lemma:compact_1}
  Let $\{v_n\}\subset S_m\cap X_2$ be any bounded Palais-Smale sequence of the constrained functional $I_{|S_m\cap X_2}$, at an arbitrary level $c>0$, satisfying $P(u_n) \to 0$. Then there exists $v \in S_m\cap X_2$ and $\mu>0$ such that, up to the extraction of a subsequence, $v_n \to v$ strongly in $H^1(\mathbb{R}^N)$ and $-\Delta v+\mu v =f(v)$.
\end{lemma}

Arguing as the proof of Lemma \ref{lemma:unbounded}, one can also establish a ``nonradial" variant in $X_2$. Using that version and repeating the argument of Lemma \ref{lemma:Emk}, we obtain
  \begin{lemma}\label{lemma:Emk_1}
     $\overline{E}_{m,k}\to+\infty$ as $k\to\infty$.
  \end{lemma}

\smallskip
\noindent
\textbf{End of the proof of Theorem \ref{theorem:nonradial2}.} For each $k\in\mathbb{N}^+$, by Lemmas \ref{lemma:J_3} and \ref{lemma:GkEmk_1}, one can find a Palais-Smale sequence $\{v^k_n\}^\infty_{n=1}\subset \mathcal{P}_m\cap X_2$ of the constrained functional $I_{|S_m\cap X_2}$ at the level $\overline{E}_{m,k}>0$. By Lemma \ref{lemma:I_4} $(iv)$, this sequence is bounded in $X_2$ and thus in view of Lemma \ref{lemma:compact_1}, we deduce that \eqref{P_m} has a nonradial solution $v_k\in X_2$ with $I(v_k)=\overline{E}_{m,k}$. Also, from Lemma \ref{lemma:GkEmk_1} $(ii)$ and Lemma \ref{lemma:Emk_1}, it follows that
  \begin{linenomath*}
    \begin{equation*}
      I(v_{k+1})\geq I(v_k)>0\qquad\text{for any}~k\geq1
    \end{equation*}
  \end{linenomath*}
and $I(v_k)\to+\infty$. \hfill $\square$

\subsection{Proof of Theorem \ref{theorem:nonradial1}}\label{subsect:nonradial1}
Recall that $N\geq4$, $X_1:=H^1_{\mathcal{O}_1}\cap X_\omega$ and $f$ is odd satisfying $(f0)-(f5)$.  For any $m>0$, we define the infimum
  \begin{linenomath*}
    \begin{equation*}
      \overline{E}_m:=\inf_{v\in\mathcal{P}_m\cap X_1}I(v),
    \end{equation*}
  \end{linenomath*}
which is positive by Lemma \ref{lemma:I_4} $(iii)$ and satisfies
  \begin{lemma}\label{lemma:Em_1}
    $\overline{E}_m > 2 E_m$.
  \end{lemma}
\proof Let $v\in \mathcal{P}_m \cap X_1$ be arbitrary. We define
\begin{linenomath*}
  \begin{equation*}
    \Omega_1:=\{x\in\mathbb{R}^N~|~|x_1|>|x_2|\}\qquad\text{and}\qquad\Omega_2:=\{x\in\mathbb{R}^N~|~|x_1|<|x_2|\}.
  \end{equation*}
\end{linenomath*}
It is clear that $\chi_{\Omega_j}v\in S_{m/2}\cap H^1_0(\Omega_j) \subset S_{m/2}$, $j=1,2$. Since
\begin{linenomath*}
  \begin{equation*}
    0=P(v)=P(\chi_{\Omega_1}v)+P(\chi_{\Omega_2}v)=2P(\chi_{\Omega_1}v),
  \end{equation*}
\end{linenomath*}
we have $\chi_{\Omega_1}v\in\mathcal{P}_{m/2}$ and thus
\begin{linenomath*}
  \begin{equation*}
    I(v)=I(\chi_{\Omega_1}v)+I(\chi_{\Omega_2}v)=2I(\chi_{\Omega_1}v)\geq 2E_{m/2}.
  \end{equation*}
\end{linenomath*}
Since $v$ is arbitrary and the function $m \mapsto E_m$ is strictly decreasing by Theorem \ref{theorem:Em}, we obtain
  \begin{linenomath*}
    \begin{equation*}
      \overline{E}_m:=\inf_{v\in\mathcal{P}_m\cap X_1}I(v)\geq 2E_{m/2}>2E_m.
    \end{equation*}
  \end{linenomath*}
The proof of the lemma is complete. \hfill $\square$

Note that, for any solution $w\in X_1$ of \eqref{P_m}, one has $w\in \mathcal{P}_m\cap X_1$ and thus $I(w)\geq \overline{E}_m$. To complete the proof of Theorem \ref{theorem:nonradial1}, it only remains to show that $\overline{E}_m$ is reached by some solution $v\in X_1$ of \eqref{P_m}. When $N-2M=0$, we have $X_1=X_2$ (with $N-2M\neq 1$). Since in that case $\overline{E}_m$ coincides with the minimax value $\overline{E}_{m,1}$ defined by \eqref{eq:Emk_1},  the result follows from the fact, shown in Subsection \ref{subsect:nonradial2}, that $\overline{E}_{m,1}$ is reached by a solution $v_1\in X_2$ of \eqref{P_m}. The rest of the proof is devoted to deal with the case $N-2M\neq 0$.

First note that adapting the proof of Lemma \ref{lemma:J_1} we can derive the following ``nonradial" version.

\begin{lemma}\label{lemma:J_4}
  Let $\mathcal{G}$ be a homotopy stable family of compact subsets of $S_m\cap X_1$ (with $B=\emptyset$) and set
    \begin{linenomath*}
      \begin{equation*}
        E_{m,\mathcal{G}}:=\inf_{A\in\mathcal{G}}\max_{v\in A}J(v).
      \end{equation*}
    \end{linenomath*}
  If $E_{m,\mathcal{G}}>0$, then there exists a Palais-Smale sequence $\{v_n\}\subset \mathcal{P}_m\cap X_1$ for the constrained functional $I_{|S_m\cap X_1}$ at the level $E_{m,\mathcal{G}}$.
\end{lemma}

Let $\overline{\mathcal{G}}$ be the class of all singletons included in $S_m\cap X_1$. Clearly, it is a homotopy stable family of compact subsets of $S_m\cap X_1$ (with $B=\emptyset$) and
  \begin{linenomath*}
    \begin{equation*}
      E_{m,\overline{\mathcal{G}}}=\inf_{v\in S_m\cap X_1}I(s(v)\star v)=\overline{E}_m>0.
    \end{equation*}
  \end{linenomath*}
Applying Lemma \ref{lemma:J_4} to $\overline{\mathcal{G}}$, we obtain
\begin{lemma}\label{lemma:I_6}
  There exists a Palais-Smale sequence $\{v_n\}\subset \mathcal{P}_m\cap X_1$ for the constrained functional $I_{|S_m\cap X_1}$ at the level $\overline{E}_m$.
\end{lemma}

To study the convergence of the Palais-Smale sequence guaranteed by Lemma \ref{lemma:I_6}, we need the following Lions type result whose proof can be found, for example, in \cite[Corollary 3.2]{Me17}.
  \begin{lemma}\label{lemma:lions}
    Assume that $N\geq4$ and $N-2M\neq0$. Let $\{v_n\}$ be a bounded sequence in $H^1_{\mathcal{O}_1}(\mathbb{R}^N)$ which satisfies, for all $r >0$,
      \begin{linenomath*}
        \begin{equation*}
          \underset{n\to\infty}{\lim}\underset{y\in\{0\}\times\{0\}\times\mathbb{R}^{N-2M}}{\sup}\int_{B(y,r)}|v_n|^2dx=0.
        \end{equation*}
      \end{linenomath*}
    Then $v_n\to0$ in $L^p(\mathbb{R}^N)$ for any $2<p<2N/(N-2)$.
  \end{lemma}

We shall also use Lemma \ref{lemma:Em_2} which follows from an adaptation of the arguments of Lemmas \ref{lemma:Em_nonincreasing} and \ref{lemma:Em_strictlymonotonic}.
\begin{lemma}\label{lemma:Em_2}
  Assume that $N\geq4$, $N-2M\neq0$, and  $f$ is an odd function satisfying $(f0)-(f4)$.  Then the following statements hold.
    \begin{itemize}
      \item[$(i)$] The function $m\mapsto \overline{E}_m$ is nonincreasing on $(0, \infty)$.
      \item[$(ii)$] If there exists $v\in S_m\cap X_1$ and $\mu\in\mathbb{R}$ such that
                      \begin{linenomath*}
                        \begin{equation*}
                          -\Delta v+\mu v=f(v)
                        \end{equation*}
                      \end{linenomath*}
      and $I(v)=\overline{E}_m$, then $\mu\geq0$. If in addition $\mu>0$, then $\overline{E}_m>\overline{E}_{m'}$ for any $m'>m$.
    \end{itemize}
\end{lemma}

With Lemmas \ref{lemma:lions} and \ref{lemma:Em_2} in hand and with the understanding that a bounded sequence $\{v_n\}\subset X_1$ is vanishing if, for all $r>0$,
\begin{linenomath*}
    \begin{equation*}
      \underset{n\to\infty}{\limsup}\underset{y\in\{0\}\times\{0\}\times\mathbb{R}^{N-2M}}{\sup}\int_{B(y,r)}|v_n|^2dx=0,
    \end{equation*}
  \end{linenomath*}
modifying accordingly the proof of Lemma \ref{lemma:decomposition}, we have the following compactness result.
\begin{lemma}\label{lemma:decomposition_1}
  Let $\{v_n\}\subset S_m \cap X_1$ be any bounded Palais-Smale sequence of the constrained functional $I_{|S_m\cap X_1}$, at the level $\overline{E}_m$, satisfying $P(u_n) \to 0$. Then there exists $v \in S_m\cap X_1$ and $\mu>0$ such that, up to the extraction of a subsequence and up to translations in $\{0\}\times\{0\}\times\mathbb{R}^{N-2M}$, $v_n \to v$ strongly in $H^1(\mathbb{R}^N)$ and $-\Delta v+\mu v =f(v)$.
\end{lemma}

\smallskip
\noindent
\textbf{End of the proof of Theorem \ref{theorem:nonradial1}.} When $N-2M\neq0$, by Lemma \ref{lemma:I_6}, we have a Palais-Smale sequence $\{v_n\}\subset \mathcal{P}_m\cap X_1$ for the constrained functional $I_{|S_m\cap X_1}$ at the level $\overline{E}_m$. By Lemma \ref{lemma:I_4} $(iv)$, this sequence is bounded in $X_1$ and thus applying Lemma \ref{lemma:decomposition_1}, we see that $\overline{E}_m$ is reached by a solution $v\in X_1$ of \eqref{P_m}. At this point, the proof of Theorem \ref{theorem:nonradial1} is complete.  \hfill $\square$

\section{Final remarks}\label{sect:final}

In this last section we justify Remark \ref{remark:Em} $(ii)$  and present two open problems.

\begin{remark}\label{remark:mu_zero}
  It has been proved in Theorem \ref{theorem:groundstate} $(ii)$ that, when $f$ is odd and $N=3,4$, we can obtain a positive ground state without assuming the condition $(f5)$. Our argument there, see Lemma \ref{lemma:decomposition}, relies on the use of a Liouville type result which allows to show that for a suspected nonnegative ground state the Lagrange multiplier is strictly positive. We shall present here an example which shows that, when $N \geq 5$ and $f$ is odd only satisfying $(f0)-(f4)$,  there exist positive ground states associated to the null Lagrange multiplier. Indirectly, this example demonstrates that, to prove the existence of ground states, the strategy developed in our paper fails for general nonlinearities when $N \geq 5$ and $(f5)$ does not hold.  It is thus an open problem to figure out if an alternative approach, not relying on the sign of the Lagrange multiplier, would give more general existence results.

  We now construct the example. For $N\geq3$, let
    \begin{linenomath*}
      \begin{equation*}
        U(x):=\frac{\big[N(N-2)\big]^{\frac{N-2}{4}}}{\big(1+|x|^2\big)^\frac{N-2}{2}},\qquad\qquad p_N:=2+\frac{4}{N}+\frac{8}{N^2}
      \end{equation*}
    \end{linenomath*}
  and $p\in (p_N,2^*)$. We define the odd continuous function
    \begin{linenomath*}
      \begin{equation*}
        f(t):=
              \left\{
                \begin{aligned}
                  &|t|^{2^*-2}t,&\qquad\text{for}~|t|\leq1,\\
                  & |t|^{p-2}t,&\qquad\text{for}~|t|>1,
                \end{aligned}
              \right.
      \end{equation*}
    \end{linenomath*}
  which satisfies $(f0)-(f4)$ but not $(f5)$. When $N\geq5$, we have $U \in L^2(\R^N)$ and for any $m \geq m_N:=N(N-2)\|U\|^2_{L^2(\mathbb{R}^N)}$ there exists a unique $\varepsilon = \varepsilon(m) > 0$ such that
	$$
      U_\varepsilon(x):= \varepsilon^{(2-N)/4}U\Big(\frac{x}{\sqrt{\varepsilon}}\Big)=\frac{\big[N(N-2)\varepsilon\big]^{\frac{N-2}{4}}}{\big(\varepsilon+|x|^2\big)^\frac{N-2}{2}}
    $$
  satisfies $\|U_\varepsilon\|^2_{L^2(\mathbb{R}^N)}=\varepsilon\|U\|^2_{L^2(\mathbb{R}^N)}=m$ and
    $$
      0<U_\varepsilon (x)\leq 1\qquad\text{for all}~x\in\mathbb{R}^N.
    $$
  In addition, it can be checked, see for example \cite[Lemma A.8]{G93}, that
	$$-\Delta U_\varepsilon = |U_\varepsilon|^{2^*-2}U_\varepsilon \qquad \text{for all}~x \in \mathbb{R}^N.$$
  Thus, in this case, \eqref{P_m} has a positive radial solution $U_{\varepsilon}$ with $\mu=0$.  We next show that $U_{\varepsilon}$ is a ground state. Denoted by $\mathcal{S}$ the best Sobolev constant such that
    \begin{equation}\label{eq:Sobolev}
      \mathcal{S} \|u\|^{2}_{L^{2^*}(\mathbb{R}^N)} \leq \|\nabla u\|^2_{L^2(\mathbb{R}^N)}\qquad\text{for any}~u\in \mathcal{D}^{1,2}(\mathbb{R}^N).
    \end{equation}
  One may note that
    \begin{linenomath*}
      \begin{equation*}
        I(U_\varepsilon)=\frac{1}{N}\|\nabla U_{\varepsilon}\|^2_{L^2(\mathbb{R}^N)} = \frac{1}{N}\mathcal{S}^{\frac{N}{2}},
      \end{equation*}
    \end{linenomath*}
  see \cite[Lemma A.8]{G93} for the second equality. Clearly, we now only need to prove that
      \begin{equation*}
          E_m:=  \inf_{u \in P_m}I(u)\geq \frac{1}{N}\mathcal{S}^{\frac{N}{2}}.
      \end{equation*}
  In view of the fact that
      \begin{equation*}
        F(t)\leq \frac{1}{2^*}|t|^{2^*}\qquad \text{for any}~t\in\mathbb{R}
      \end{equation*}
  and \eqref{eq:Sobolev}, it is not difficult to deduce that
      \begin{equation*}
        \begin{aligned}
          \inf_{u \in S_m}\max_{s \in \mathbb{R}}I(s \star u)
          &\geq \inf_{u \in S_m}\max_{s \in \mathbb{R}}\left[\frac{1}{2}\int_{\mathbb{R}^N}|\nabla (s \star u)|^2dx-\frac{1}{2^*}\int_{\mathbb{R}^N}|s\star u|^{2^*}dx\right]\\
          &=\inf_{u \in S_m}\max_{s \in \mathbb{R}}\left[\frac{1}{2}e^{2s}\int_{\mathbb{R}^N}|\nabla u|^2dx-\frac{1}{2^*}e^{2^*s}\int_{\mathbb{R}^N}|u|^{2^*}dx\right]\\
          &=\inf_{u \in S_m}\left[\frac{1}{2}\left(\frac{\|\nabla u\|^2_{L^2}}{\|u\|^{2^*}_{L^{2^*}}}\right)^{\frac{2}{2^*-2}}\|\nabla u\|^2_{L^2}-\frac{1}{2^*}\left(\frac{\|\nabla u\|^2_{L^2}}{\|u\|^{2^*}_{L^{2^*}}}\right)^{\frac{2^*}{2^*-2}}\|u\|^{2^*}_{L^{2^*}}\right]\\
          &= \inf_{u \in S_m}\frac{1}{N}\left(\frac{\|\nabla u\|^2_{L^2}}{\|u\|^2_{L^{2^*}}}\right)^{\frac{N}{2}}\geq \frac{1}{N}\mathcal{S}^{\frac{N}{2}}.
        \end{aligned}
      \end{equation*}
	Combining Lemma \ref{lemma:I_3}, we obtain
	  \begin{linenomath*}
        \begin{equation*}
         E_m:=  \inf_{u \in P_m}I(u)=\inf_{u \in S_m}\max_{s \in \mathbb{R}}I(s \star u)\geq \frac{1}{N}\mathcal{S}^{\frac{N}{2}}
        \end{equation*}
      \end{linenomath*}
    and thus $U_{\varepsilon}$ is a ground state. This example shows that, when $N \geq 5$ and for an arbitrary nonlinearity satisfying $(f0)-(f4)$, the proof of Lemma \ref{lemma:decomposition} breaks down since there is not hope to show that the Lagrange multiplier,  whose value is given in \eqref{eq:mu}, is strictly positive.
\end{remark}

\begin{remark}\label{remark:E_infinity}
  When $N\geq3$ and $f$ satisfies $(f0)-(f4)$, Lemmas \ref{lemma:I_4} and \ref{lemma:Em_nonincreasing} imply that
    \begin{linenomath*}
      \begin{equation*}
        E_\infty := \lim_{m \to \infty} E_m
      \end{equation*}
    \end{linenomath*}
  exists and $E_\infty \geq 0$. In particular, if $(f6)$ also holds, then $E_\infty=0$ by Lemma \ref{lemma:Em_atinfinity}. Let us show that, when $(f6)$ is replaced by the somehow opposite condition
    \begin{itemize}
      \item[$(f6)'$] $\limsup_{t \to 0} f(t)t/|t|^{\frac{2N}{N-2}}<+\infty$,
    \end{itemize}
  then $E_{\infty} >0$. Indeed, under the conditions $(f0)-(f4)$ and $(f6)'$, it is possible to define $I$ as a free functional on $\mathcal{D}^{1,2}(\mathbb{R}^N)$. In addition $I$ is of class $C^1$ and setting
    \begin{linenomath*}
      \begin{equation*}
        c_{mp}:=\inf_{\gamma \in \Gamma}\max_{t \in [0,1]}I(\gamma(t)),
      \end{equation*}
    \end{linenomath*}
  where $\Gamma:=\{\gamma \in C([0,1],\mathcal{D}^{1,2}(\mathbb{R}^N))~|~\gamma(0)=0,~I(\gamma(t))<0\}$  one has that $c_{mp}>0$. Clearly, if we show that
    \begin{linenomath*}
      \begin{equation}\label{eq:remark_1}
        E_m \geq c_{mp}\qquad\text{for any}~m>0,
      \end{equation}
    \end{linenomath*}
  then the proof is complete. In order to get \eqref{eq:remark_1}, we adapt the argument of \cite[Lemma 7.1]{BJL13}. For any given $u \in \mathcal{P}_m$, in view of Lemma \ref{lemma:I_1} $(i)$, there exists $s_-<0$ such that
    \begin{linenomath*}
      \begin{equation}\label{eq:remark_2}
        I(\theta (s_-\star u))\leq e^{2s_-}\int_{\mathbb{R}^N}|\nabla u|^2dx<I(u)\qquad\text{for any}~\theta\in [0,1].
      \end{equation}
    \end{linenomath*}
  By Lemma \ref{lemma:I_2} $(ii)$, we can also choose $s_+>0$ large enough such that $I(s_+ \star u)<0$. Since the path
    \begin{linenomath*}
      \begin{equation*}
        \gamma(t):=
          \left\{
              \begin{aligned}
                2t (s_-\star u),\qquad\qquad&\qquad 0\leq t \leq \frac{1}{2},\\
                [2(1-t)s_-+(2t-1)s_+]\star u,&\qquad \frac{1}{2}\leq t \leq 1,
              \end{aligned}
          \right.
      \end{equation*}
    \end{linenomath*}
  belongs to $\Gamma$, from \eqref{eq:remark_2} and Lemma \ref{lemma:I_3} $(ii)$, it follows that
    \begin{linenomath*}
      \begin{equation*}
        I(u)=\max_{t\in [0,1]}I(\gamma(t))\geq c_{mp}>0.
      \end{equation*}
    \end{linenomath*}
  Noting that $u \in \mathcal{P}_m$ is arbitrary, we obtain \eqref{eq:remark_1} and thus $E_\infty \geq c_{mp}>0$.

  As an example of function that satisfies $(f0)-(f5)$ and $(f6)'$, we have
    \begin{linenomath*}
      \begin{equation*}
        f(t):=\beta\left[1-\frac{\beta_N(N-2)|t|^{\beta_N}}{2N(1+|t|^{\beta_N})} \right] \frac{|t|^{\frac{4}{N-2}}t}{1+|t|^{\beta_N}}
      \end{equation*}
    \end{linenomath*}
  with its primitive integral
    \begin{linenomath*}
      \begin{equation*}
        F(t):=\beta\frac{(N-2)|t|^{\frac{2N}{N-2}}}{2N(1+|t|^{\beta_N})},
      \end{equation*}
    \end{linenomath*}
  where $\beta >0$ and $\beta_N \in (0, \frac{4}{N(N-2)}]$.
\end{remark}

\begin{remark}\label{remak:conjecture}
  To remind the dependence of $E_m$ and $E_\infty$ on $f$, we now denote them as $E_{f,m}$ and
    \begin{linenomath*}
      \begin{equation*}
        E_{f,\infty}:=\lim_{m\to\infty}E_{f,m}.
      \end{equation*}
    \end{linenomath*}
  Let $N\geq3$ and consider the functions that satisfy $(f0)-(f5)$ and
    \begin{itemize}
      \item[$(f6)''$] $\lim_{t\to0}f(t)t/|t|^{\frac{2N}{N-2}}=:L_f\in (0,+\infty]$.
    \end{itemize}
  In view of Remark \ref{remark:E_infinity}, for any $f,g$, we conjecture that
    \begin{linenomath*}
      \begin{equation*}
        E_{f,\infty}>E_{g,\infty}~\text{if}~L_f<L_g,\qquad~E_{f,\infty}=E_{g,\infty}~\text{when}~L_f=L_g;
      \end{equation*}
    \end{linenomath*}
  or, at least, $E_{f,\infty}\geq E_{g,\infty}$ if $L_f<L_g$. Clearly, by Remark \ref{remark:E_infinity}, this conjecture is true when $L_g=+\infty$.
\end{remark}


\section*{Acknowledgment}
\addcontentsline{toc}{section}{Acknowledgment}

The authors thank  Jaroslaw Mederski for pointing to them the reference \cite{BM20}.  This has led  us to improve a first version of our paper and, in particular,  to show that the ground state obtained in Theorem \ref{theorem:groundstate} $(ii)$ can be assumed to be radially symmetric.  S.-S. Lu acknowledges the support of the National Natural Science Foundation of China (NSFC-11771324, 11831009 and 11811540026), of the China Scholarship Council (CSC-201706250149) and the hospitality of the Laboratoire de Math\'{e}matiques (CNRS UMR 6623), Universit\'{e} de Bourgogne Franche-Comt\'{e}.

\end{document}